\documentclass[10pt]{elsart_CAMCS}

\usepackage{savesym}
\savesymbol{AND}
\usepackage{amsmath,amssymb,mathrsfs,array,stmaryrd,natbib,color,algorithmic,algorithm}
\usepackage[xdvi]{graphicx}

\newcommand{\totpolspace}{{\mathbb{P}_{\No}}}
\newcommand{\Nd}	{d}
\newcommand{\Ndx}	{{\Nd_{\boldsymbol{x}}}}
\newcommand{\NdF}	{{d_F}}
\newcommand{\Ndnu}	{{d_\nu}}
\newcommand{\Nq}	{{N_q}}
\newcommand{\Nqt}	{\widetilde{N_q}}
\newcommand{\Nqhat}	{\widehat{N_q}}
\newcommand{\bx}	{\boldsymbol{x}}
\newcommand{\xq}	{x^{(q)}}
\newcommand{\bxq}	{\bx^{(q)}}
\newcommand{\xqt}	{x^{(\widetilde{q})}}
\newcommand{\bxi}	{\boldsymbol{\xi}}
\newcommand{\bzeta}	{\boldsymbol{\zeta}}
\newcommand{\bzetaq}	{\bzeta^{(q)}}

\newcommand{\xiq}	{\xi^{(q)}}

\newcommand{\uq}	{u^{(q)}}
\newcommand{\bu}	{\boldsymbol{u}}

\newcommand{\uhat}	{\widehat{u}}
\newcommand{\buhat}	{\widehat{\bu}}
\newcommand{\uqt}	{{u^{(\widetilde{q})}}}
\newcommand{\uqh}	{{u^{(\widehat{q})}}}

\newcommand{\wq}	{{w^{(q)}}}
\newcommand{\bxiq}	{{\bxi^{(q)}}}
\newcommand{\bxiqt}	{{\bxi^{(\widetilde{q})}}}
\newcommand{\bxiqh}	{{\bxi^{(\widehat{q})}}}
\newcommand{\matPsi}	{\Psi}

\newcommand{\coef}	{c}

\newcommand{\bcoef}	{\boldsymbol{c}}
\newcommand{\bcoeft}	{\widetilde{\boldsymbol{c}}}

\newcommand{\No}	{p}
\newcommand{\NoLARS}	{\widetilde{\No}}
\newcommand{\Not}	{\widetilde{p}}
\newcommand{\balpha}	{{\boldsymbol{\alpha}}}

\newcommand{\bgamma}	{{\boldsymbol{\gamma}}}

\newcommand{\nr}	{{n_r}}

\newcommand{\Ninter}	{{N_l}}			
\newcommand{\NinterPC}	{{N_l^{\rm (PC)}}}			

\newcommand{\group}	{\mathcal{J}}

\newcommand{\modex}	{w}
\newcommand{\modesto}	{\lambda}
\newcommand{\bmodex}	{{\boldsymbol{\modex}}}
\newcommand{\bmodesto}	{\boldsymbol{\modesto}}
\newcommand{\res}	{z}
\newcommand{\resq}	{\res^{(q)}}
\newcommand{\bres}	{\boldsymbol{\res}}

\newcommand{\nJ}	{n_f}

\newcommand{\nJzero}	{n_{f_1}}

\newcommand{\borel}	{\mathcal{B}}
\newcommand{\borelxi}	{\mathcal{B}_{\Xi}}
\newcommand{\measbrut}	{\mu}
\newcommand{\meas}	{\mu_\Theta}
\newcommand{\measxi}	{{\mu_{\Xi}}}
\newcommand{\measxipasbold} {{\mu_{\Xi}}}
\newcommand{\measxioneD}{{\mu_{\Xi}}}

\newcommand{\randomspace} {\mathcal{S}}
\newcommand{\randomspaceNo} {\randomspace_\No}

\newcommand{\basisset}	{\mathcal{J}}
\newcommand{\basissetf}	{\mathcal{J}_f}

\newcommand{\basissetprior}	{\basisset_{\rm prior}}

\newcommand{\basissetpost}	{\basisset_{\rm post}}
\newcommand{\basissetfpost}	{\basisset_{f, {\rm post}}}
\newcommand{\basissetcur}	{\basisset_{\rm eff}}
\newcommand{\basissetfcur}	{\basisset_{f, {\rm eff}}}

\newcommand{\cardP}	{{|\basisset|}}
\newcommand{\cardPpost}	{{|\basisset_{\rm post}|}}

\newcommand{\cardPprior}	{{|\basisset_{\rm prior}|}}

\newcommand{\cardPcur}	{{|\basisset_{\rm eff}|}}
\newcommand{\cardPfcur}	{{|\basisset_{f, {\rm eff}}|}}

\newcommand{\cardPgroup}	{{|\group_{\bgamma}|}}			
\newcommand{\cardPgroupmean}	{{\overline{|\group_{\bgamma, {\rm post}}|}}}		
\newcommand{\cardx}	{{|\basisset_x|}}
\renewcommand{\span}	{{\rm{span}}}
\newcommand{\proj}	{\mathcal{P}}
\newcommand{\ft}	{\widetilde{f}}
\newcommand{\fhat}	{\widehat{f}}
\newcommand{\vx}	{v_{1}}
\newcommand{\vy}	{v_{2}}
\newcommand{\Su}	{S_{\vx}}
\newcommand{\Sv}	{S_{\vy}}
\newcommand{\Sh}	{{S_h}}
\newcommand{\NH}	{{N_a}}
\newcommand{\KLrank}	{N}

\newcommand{\apriori}	{{\textit{a priori}}}
\newcommand{\aposteriori}{{\textit{a posteriori}}}
      \newcommand{\ieLM}	{\textit{i.e.}}
      \newcommand{\egLM}	{\textit{e.g.}}
      
\newcommand{\MC}	{{Monte Carlo}}

\renewcommand{\d}	{{\rm{d}}}

\newcommand{\be}	{\begin{equation}}
\newcommand{\ee}	{\end{equation}}
\newcommand{\bi}	{\begin{itemize}}
\newcommand{\ei}	{\end{itemize}}
\newcommand{\bea}	{\begin{eqnarray}}
\newcommand{\eea}	{\end{eqnarray}}
\newcommand{\bc}	{\begin{center}}
\newcommand{\ec}	{\end{center}}

\DeclareMathOperator*{\argmin}	{{\rm arg\,min}}

\DeclareMathOperator*{\esp}	{\mathscr{E}}
\DeclareMathOperator*{\Var}	{{\rm Var}}
\DeclareMathOperator*{\Ident}	{I}	
\DeclareMathOperator*{\vectorize}{\mathrm{vec}}	
\newcommand{\lra}	{\longrightarrow}

\newtheorem{remark}{Remark}

\begin{document}

\begin{frontmatter}
\title {Quantification of uncertainty from high-dimensional scattered data via polynomial approximation}

\author[label1,label2]{Lionel Mathelin}		
\ead{mathelin@limsi.fr}
\address[label1]{LIMSI-CNRS, BP 133, 91403 Orsay, France.}
\address[label2]{Dpt. of Aeronautics and Astronautics, Massachusetts Institute of Technology, \\ 77 Massachusetts Av., Cambridge, MA 02139, USA.}


\begin{abstract}
This paper discusses a methodology for determining a functional representation of a random process from a collection of scattered pointwise samples.
The present work specifically focuses onto random quantities lying in a high dimensional stochastic space in the context of limited amount of information. The proposed approach involves a procedure for the selection of an approximation basis and the evaluation of the associated coefficients. The selection of the approximation basis relies on the {\apriori} choice of the High-Dimensional Model Representation format combined with a modified Least Angle Regression technique. The resulting basis then provides the structure for the actual approximation basis, possibly using different functions, more parsimonious and nonlinear in its coefficients. To evaluate the coefficients, both an alternate least squares and an alternate weighted total least squares methods are employed. Examples are provided for the approximation of a random variable in a high-dimensional space as well as the estimation of a random field. Stochastic dimensions up to 100 are considered, with an amount of information as low as about 3 samples 
per dimension, and robustness of the approximation is demonstrated w.r.t. noise in the dataset. The computational cost of the solution method is shown to scale only linearly with the cardinality of the {\apriori} basis and exhibits a $\left(\Nq\right)^s$, $2 \le s \le 3$, dependence with the number $\Nq$ of samples in the dataset. The provided numerical experiments illustrate the ability of the present approach to derive an accurate approximation from scarce scattered data even in the presence of noise.
\end{abstract}


\end{frontmatter}




\section{Introduction}

With the growing available computational power, and as more efficient numerical methods become available, domains as diverse as engineering, chemistry, psychometrics, medicine, finance or social sciences, now heavily rely on simulation for the prediction of more and more complex phenomena, often combining multi-models and high accuracy requirement. The prediction capability of modern simulations is often such that a new bottleneck for accuracy has emerged from the lack of relevant boundary and/or initial conditions (BICs) as well as parameters intrinsic to the model of the system at hand, \egLM, diffusivity, viscosity, etc. These sources of uncertainty are hereafter simply referred to as BICs. They are often poorly known and have to be estimated or modeled. This introduces modeling errors which often constitute the main source of lack of accuracy in the simulation chain. This situation has triggered a renewed interest for stochastic modeling where it is explicitly accounted for uncertainty in the model.
The BICs may sometimes be modeled from first principles but are often approximated in a functional form involving a set of influencing parameters and identified from experimental measurements. However, more often than not, only relatively few measurements are available, in particular when a significant number of parameters is of influence so that representing the BICs takes the form of a high-dimensional approximation problem.

If the random process, which output is to be represented in closed-form, is driven by known equations, efficient techniques may be used to determine its representation. In the specific case of high-dimensional quantities, tensor-based representations have proved to be effective when applicable. In particular, low-rank approximations based on an {\apriori} chosen separated representation can be efficiently derived, see \cite{Nouy_CMAME_07, Nouy_CMAME_10, Nouy_HD_10, Matthies_Zander_12} in the context of uncertainty quantification (UQ).
If a closed-form model description of the process at hand is not available, one is typically left with approximating it from a finite collection of instances, hereafter termed samples. When the process is known only from a closed numerical code used as a black-box or if measurements can be made arbitrarily (design of experiments), some properties of approximation theory can be exploited. For instance, measurements may be taken at some particular locations in the parameter space, possibly associating a weight to them, so that the random Quantity of Interest (QoI) can be represented in the retained approximation basis with good accuracy using (sparse) quadrature techniques, \cite{Novak_Ritter_99}; see also \cite{Xiu_Hesthaven_05} for an application to UQ. Anisotropy in the QoI may be exploited by biasing the quadrature weights, \cite{Nobile_al_07, Ganapathysubramanian_Zabaras_07, Ma_Zabaras_10}. In \cite{Doostan_Iaccarino_09}, an Alternate Least Squares (ALS) technique to estimate the coefficients has been 
considered with samples lying on a tensor-product grid. Another situation of design of experiment arises in importance sampling where the Markov-Chain {\MC} algorithm requires a new sample at a specific proposed location. This control over the samples usually brings efficiency and allows to approximate a reasonably behaved QoI with accuracy.

A different situation occurs when the data are scattered, with no ability to choose the set of samples nor to add a measurement. This is a common situation, typically arising when samples come from a past experiment or are costly to acquire so that new samples cannot be taken. In this context, one has to resort to a regression-based approach and the coefficients of the approximation are then solution of an optimization problem. This type of approach was considered in \cite{Choi_etal_2004, Berveiller_etal_2006, Beylkin_etal_2009}.

In the present work, the focus is specifically put on deriving a closed-form approximation of a high-dimensional quantity of interest from a small, uncontrolled, collection of its samples.
This requires to determine an approximation basis finely tuned to the data at hand and an efficient way of evaluating the associated coefficients.
To this aim, we rely on the fact that, as a counterpart of the curse of dimensionality associated with high-dimensional problems, real applications often reward with a \emph{blessing} of dimensionality. Indeed, in many cases, the QoI can be well approximated in a low-dimensional subspace of the solution space, sometimes involving orders of magnitude fewer degrees-of-freedom. This typically occurs when the solution exhibits some degree of sparsity in the retained functional space. Efficient techniques have been proposed in the recent past to take advantage of this situation and essentially consist in matching the approximation with the observational data while promoting a sparse coefficient set. This class of methods work well in many different contexts and have been recently applied to the UQ framework, \cite{Doostan_Owhadi_11, Mathelin_Gallivan_12}. These techniques rely on the Compressed Sensing theory, \egLM, \cite{Candes_Tao_decoding,Donoho_06}, and may seem well suited for the present problem as they 
promote a low cardinality approximation of the QoI. However, they require to handle a potentially huge representation basis, or \emph{dictionary}, and associated optimization problem, leading to severe memory and computation limitations in the present high-dimensional context.

In this paper, we present a solution method combining the strength of different techniques, taking advantage of the sparsity of the representation in a suitable basis and allowing an efficient approximation of a well-behaved multivariate function with a low number of degrees-of-freedom hence compatible with a small experimental dataset. The driving principle is first to consider a tight approximation basis based on {\apriori} knowledge on the QoI at hand and to rely on the available data to further refine it.
In a nutshell, an initial approximation basis is first considered in the High-Dimensional Model Representation format (HDMR, \cite{Rabitz_Alis_99,Alis_Rabitz_01}), assuming it is suitable for representing the QoI. This initial basis is hereafter referred to as {\apriori} basis. Next, available data are used to refine it by retaining only its most relevant basis functions through a constructive subset selection procedure based on a modification of the Least Angle Regression approach proposed in \cite{Efron_etal_04}. This {\aposteriori} basis defines a skeleton from which a final basis is built and the associated coefficients are evaluated with an alternate least squares technique. The solution method allows to approximate random variables as well as random fields and is here shown to outperform both sparse grids and tensored-based techniques.

The paper is organized as follows. The representation of a random quantity is central to the methodology discussed in this paper. Standard techniques for deriving a closed-form approximation of a random variable from a finite set of samples are briefly recalled in section \ref{UQ_section}. Similarly, different representation formats of functions in high-dimensional spaces are subsequently heavily used in the paper and a short discussion is given in section \ref{representation_section}. The proposed solution method is introduced and discussed in section \ref{methodology_section} and an algorithm is given. Scalability of the proposed approach together with its robustness w.r.t. noise in the data is also discussed. In section \ref{Sec_results}, the present methodology is illustrated on a stochastic diffusion equation involving up to 100 dimensions and on the space-dependent solution of the Shallow Water Equations with random parameters. Accuracy, robustness and scalability of the proposed approach are shown. 
Concluding remarks close the paper in section \ref{Conclusion_section}.

\section{Quantification of uncertainty}  \label{UQ_section}

Thanks to its pivotal role in the rest of the paper, the representation of a random quantity and standard ways of evaluating it in closed-form from a discrete set of samples is now briefly discussed.

\subsection{General framework}

Random quantities are defined on a probability space $\left(\Theta, \borel_\Theta, \meas \right)$ where $\Theta$ is the space of elementary events $\theta \in \Theta$, $\borel_\Theta$ a $\sigma$-algebra defined on $\Theta$ and $\meas$ a probability measure on $\borel_\Theta$. To make the description of the problem amenable to a tractable representation, it is convenient to introduce a finite set of statistically independent random variables $\left\{\xi_i \right\}_{i=1}^{\Nd}: \Theta \rightarrow \Xi_i \subseteq \mathbb{R}, \: \theta \mapsto \xi_i(\theta)$. The set of these $\Nd$ random variables is defined on a probability space $\left(\Xi, \borelxi, \measxi\right)$ with $\Xi = \times_{i=1}^\Nd \Xi_i = \bxi\left(\Theta\right) \subseteq \mathbb{R}^{\Nd}$, $\bxi := \left(\xi_1 \, \ldots \, \xi_\Nd\right)$, $\borelxi \subset 2^{\Xi}$ a $\sigma$-algebra on $\Xi$ and $\measxi = \meas \circ \bxi^{-1}$ the probability measure on $\borelxi$.
Since the physical process at hand relies on random quantities belonging to $\left(\Theta, \borel_\Theta, \meas \right)$, a suitable description of its output, or its solution in case the physical process is described by a known mathematical model, may be determined in $\left(\Xi, \borelxi, \measxi\right)$ as justified by the Doob-Dynkin lemma.

In this work, we restrict ourselves to random variables of physical significance, \ieLM, real-valued second order variables satisfying:
\be
\esp_\theta \left[u\left(\theta\right)^2\right] := \int_\Theta {u\left(\theta\right)^2 \, \d\meas\left(\theta\right)} = \int_\Xi {u\left(\bzeta\right)^2 \, \d\measxi\left(\bzeta\right)} =: \esp_{\bxi} \left[u\left(\bxi\right)^2\right] \qquad < + \infty,
\ee
where $\esp$ denotes the expectation operator and $u$ is the quantity of interest (QoI). It is then natural to consider the space of square integrable functions $\randomspace$ for describing real-valued functions of the random quantities:
\be
\randomspace := L^2\left(\Xi, \, \measxi\right) = \left\{ v: \Xi \rightarrow \mathbb{R}, \bxi \mapsto v\left(\bxi\right); \: \esp_{\bxi} \left[v\left(\bxi\right)^2\right] < + \infty\right\}.
\ee

Upon introduction of a natural inner product of $\randomspace$: $\displaystyle \left< v, w \right>_{L^2\left(\Xi,\, \measxi\right)} := \int_\Xi {v\left(\bzeta\right) \, w\left(\bzeta\right) \, \d\measxi\left(\bzeta\right)}$, $\forall \, v, w \in \randomspace$, and the associated norm $\left\|v\right\|^2_{L^2\left(\Xi, \, \measxi\right)} := \left<v, v\right>_{L^2\left(\Xi, \, \measxi\right)}$,
$\randomspace$ is a Hilbert space. Further, we define $\left<v\right>_{L^2\left(\Xi, \, \measxi\right)} := \esp_{\bxi} \left[v\left(\bxi\right)\right]$. One can now rely on functional analysis results and take advantage of approximation theory techniques to characterize the output $u$. Introducing a Hilbertian basis $\left\{\psi_k\right\}_{k \in \mathbb{N}}$ of $\randomspace$, the output can then be uniquely represented as $\displaystyle u\left(\bxi\right) = \sum_\alpha {c_\alpha \, \psi_\alpha\left(\bxi\right)}$.

The basis $\left\{\psi_\alpha\right\}_{\alpha \in \mathbb{N}}$ is typically chosen orthonormal w.r.t. the inner product $\left< v, w \right>_{L^2\left(\Xi, \, \measxi\right)}$.
Orthonormality of the basis leads to $\left< \psi_\alpha, \psi_{\alpha'} \right>_{L^2\left(\Xi, \,\measxi\right)} = \delta_{\alpha \alpha'}, \: \forall \, \alpha, \alpha' \in \mathbb{N}$, with $\delta$ the Kronecker delta, and the decomposition coefficients $\left\{c_\alpha\right\}$ then express as
\be
c_\alpha = \left< u, \psi_\alpha \right>_{L^2\left(\Xi, \, \measxi\right)} = \int_\Xi {u\left(\bzeta\right) \, \psi_\alpha\left(\bzeta\right) \, \d\measxi\left(\bzeta\right)}, \qquad \forall \, \alpha \in \mathbb{N}.
\label{coef_express}
\ee

For a given representation basis $\left\{\psi_\alpha\right\}$ of $\randomspace$, the output $u\left(\bxi\right)$ is entirely characterized by the set of coefficients $\left\{c_\alpha\right\}$. For computational purpose, the infinite dimensional representation is substituted with a finite dimensional approximation relying on a subset $\basisset \subset \mathbb{N}$ of the representation basis:
\be
u\left(\bxi\right) \approx \sum_{\alpha \in \basisset} {c_\alpha \, \psi_\alpha\left(\bxi\right)}. \label{approx_repres}
\ee

\subsection{Computing a data-driven approximation} \label{Comput_coefs}
As seen above, in many situations, a closed-form model of the QoI is not available or not reliable enough to be used and one can only rely on the sole available input-output information to approximate the output $u$. The solution method then consists in using a set of outputs given some inputs, \ieLM, samples of the process. One then looks for a functional form of the map between the set of random variables $\bxiq$ and the output value $u\left(\bxiq\right) =: \uq$, $\forall \, 1 \le q \le \Nq$, where $\Nq$ is the size of the available experimental set. Approximating the output under the functional form of Eq.~\eqref{approx_repres} results in evaluating the coefficients $\left\{c_\alpha\right\}$ from $\left\{ \left(\bxiq,\uq\right) \right\}_{q=1}^{\Nq}$, $\bxiq = \left(\xiq_1 \, \ldots \, \xiq_\Nd \right)$.

\subsubsection{Direct evaluation}

If the sampling can be controlled, in the sense that samples can be drawn arbitrarily, the popular {\MC} approach can be followed and the approximation coefficients are then estimated from
\be
c_\alpha = \int_\Xi {u\left(\bzeta\right) \, \psi_\alpha\left(\bzeta\right) \, \d\measxi\left(\bzeta\right)} \approx \sum_q {u\left(\bxiq\right) \, \psi_\alpha\left(\bxiq\right)}.
\ee

{\MC}-based estimation is very robust and easy to implement but suffers from a slow $\mathcal{O}\left(\Nq^{-1/2}\right)$ asymptotic convergence rate. However, since the convergence rate does not depend on the dimensionality of the integral, this is a wise choice for very high-dimensional problems where other methods fail. Alternatively, quasi-{\MC} methods generate a low-discrepancy sequence of samples improving the convergence rate of the evaluation for moderate- to high-dimensional problems.
%

For low to moderate dimensionality problems, the $\Nd$-dimensional integral arising in Eq.~\eqref{coef_express} may be advantageously evaluated with a quadrature rule:
\be
c_\alpha = \int_\Xi {u\left(\bzeta\right) \, \psi_\alpha\left(\bzeta\right) \, \d\measxi\left(\bzeta\right)} \approx \sum_q {\wq \, u\left(\bxiq\right) \, \psi_\alpha\left(\bxiq\right)},
\ee
where $\left\{\wq\right\}$ are the weights associated with the quadrature points $\left\{\bxiq\right\}$, \cite{Abramowitz_Stegun_72}.

\subsubsection{Regression} \label{subsection_regression}

The above methods require some kind of control over the samples.
If no experimental design can be exploited, a solution method is then to reformulate the evaluation of the coefficients as a minimization problem:
\be
\bcoef = \argmin_{\bcoeft \in \mathbb{R}^\cardP}{\left\| \bu - \matPsi \, \bcoeft\right\|_2}, \label{LS}
\ee
with $\bcoef = \left(c_1 \, \ldots \, c_\cardP\right)^T$, $\bu = \left( u^{(1)} \, \ldots \, u^{(\Nq)}\right)^T$, $\matPsi \in \mathbb{R}^{\Nq \times \cardP}$, $\matPsi_{q \alpha} = \psi_\alpha\left(\bxiq\right)$ and $\cardP$ the cardinality of the approximation basis $\left\{\psi_\alpha\right\}_{\alpha \in \basisset}$. For a full column rank $\matPsi$, the solution is given by $\bcoef = \matPsi^+ \, \bu$ which is typically evaluated using the Cholesky decomposition of the symmetric positive definite matrix $\matPsi^T \, \matPsi$ or the QR decomposition of $\matPsi$. When the size of the dataset grows, this standard Least Squares (LS) problem may become computationally involved. 
The quasi-regression solution alleviates the computational burden and is given by
\be
\coef_\alpha = \boldsymbol{\psi}_\alpha^T \, \bu / \left\| \psi_\alpha \right\|^2_2, \qquad \boldsymbol{\psi}_\alpha = \left(\psi_\alpha\left(\bxi^{(1)}\right) \, \ldots \, \psi_\alpha\left(\bxi^{(\Nq)}\right) \right)^T, \qquad 1 \le \alpha \le \cardP.
\ee

Standard least squares formulation as considered in Eq.~\eqref{LS} treats all predictors $\left\{\psi_\alpha\right\}_{\alpha=1}^{\cardP}$ the same way and uses the available data to estimate all the coefficients to produce an estimate with a low bias but often a large variance.
As will be discussed in section \ref{direct_approach}, additional properties of the QoI may be exploited or 
imposed to the approximation coefficients. This class of approaches trades some increase in bias with a decrease in variance and often results in an improved accuracy. A suitable solution method then typically formulates as a penalized least squares problem:
\be
\bcoef = \argmin_{\bcoeft \in \mathbb{R}^\cardP}{\left\| \bu - \matPsi \, \bcoeft\right\|_2 + \mathscr{J}\left(\bcoeft\right)}. \label{pLS}
\ee

The properties of the penalized LS solution are driven by the choice of the function $\mathscr{J}$ which flexibility leads to a variety of solution techniques, see \cite{Hesterberg_etal_08, Hastie_etal_09}. Since we have no control over the sampling strategy, we will rely on regression to estimate the approximation coefficients. The discussion of an efficient least squares formulation in the present context is postponed to section \ref{methodology_section}.

\section{Functional representation of random variables} \label{representation_section}

\subsection{Tensored bases}

As seen above, a random quantity is conveniently approximated in a Hilbertian basis $\left\{ \psi_k \right\}$.
If the random quantity is known, or expected, to exhibit a certain degree of smoothness along the stochastic space, a suitable and popular choice is to take advantage of this smoothness using a spectral-based approximation relying on polynomials. Early efforts towards this direction are the pioneering works of \cite{Wiener_38} who used univariate Hermite polynomials $\psi_\alpha\left(\xi_i\right)$ of zero-centered, unit variance, normal random variables $\xi_i \sim \mathcal{N}\left(0, 1\right)$. These polynomials define an orthogonal basis of ${L^2\left(\Xi_i, \, \measxipasbold_i\right)}$, $\measxipasbold_i \propto e^{-\frac{1}{2} \xi_i^2}$. Tensorization of univariate Hermite polynomials $\psi$ then leads to an orthogonal basis of ${L^2\left(\Xi, \, \measxi\right)}$:
\be
\left< \psi_\alpha, \psi_{\alpha'} \right>_{L^2\left(\Xi, \, \measxi\right)} \propto \int_{\Xi}{\psi_\alpha\left(\bzeta\right) \, \psi_{\alpha'}\left(\bzeta\right) \, e^{-\frac{1}{2} (\bzeta^T \bzeta)} \, \d\bzeta} \propto \delta_{\alpha \alpha'}.
\ee

This can be extended to polynomials orthogonal with respect to different measures, \cite{Ghanem_Spanos_03, Xiu_Karniadakis_02, Soize_Ghanem_04}, and constitutes the so-called (generalized) Polynomial Chaos (PC) basis. A common practice is to consider an approximation space $\randomspaceNo$ spanned by polynomials of given maximum total degree $\No$:
\be
\randomspaceNo = \span\left(\left\{\psi_\balpha\left(\bxi\right) = \psi_{\alpha_1}\left(\xi_1\right) \ldots \psi_{\alpha_\Nd}\left(\xi_\Nd\right)\right\}; \balpha = \left(\alpha_1 \, \ldots \, \alpha_\Nd\right), \sum_{i=1}^\Nd{\alpha_i} \le \No\right), \label{PC_basis_span}
\ee
and the number of terms to be determined in the approximation \eqref{approx_repres} is then $\cardP = \left( \begin{array}{c} \Nd + \No \\ \Nd \end{array}\right)$. We adopt the convention $\psi_1 \equiv 1$. When the random quantity is not smooth enough for a low degree polynomial fit to be accurate, approximation schemes such as $h/p$-type refinement or Multi-Resolution Analysis may be applied, see \cite{LeMaitre_Knio_10}.

%

Some alternative representation formats specifically exploit the tensor-product structure of the Hilbert stochastic space $\randomspace$ and approximates a $\Nd$-variate function with a series of products of lower dimensional functions. Efficient algorithms allow to determine the approximation coefficients of the representation by solving a series of low-dimensional problems while never considering the full-dimensional problem at once. A general presentation of tensor-structured numerical methods can be found in \cite{Khoromskij_review} while application to the approximation of a high-dimensional random quantity is considered in \cite{Doostan_Iaccarino_09,Nouy_HD_10,Khoromskij_Schwab_11,Matthies_Zander_12}. For instance, a $\Nd$-variate quantity may be approximated under a CANDECOMP-PARAFAC (CP) format, \cite{Harshman_70,Carroll_Chang_70}, with a sum of rank-1 terms, the simplest form of tensored-structure format:
\be
u\left(\bxi\right) \approx \sum_{r=1}^\nr{f_{1,r}\left(\xi_1\right) \, \ldots \, f_{\Nd,r}\left(\xi_\Nd\right)}, \label{CP-like}
\ee
with $\nr$ the retained rank of the decomposition and $\left\{f_{i,r}\right\}_{i=1}^\Nd$ univariate functions.
Assuming $\No$-th order polynomials for $\left\{f_{i,r}\right\}$, the resulting cardinality of the approximation is $\Nd \, \nr \, \No$. It thus exhibits a linear dependence with the number of dimensions, in contrast with the exponential dependence of the Polynomial Chaos.
Alternative decomposition techniques, easier to evaluate and numerically more stable than decomposition~\eqref{CP-like}, such as the Tucker or Tensor-Trains, can be considered, see \cite{Khoromskij_review}.
A tensored-structure format then constitutes a method of choice for deriving memory- and CPU-efficient approximation of high-dimensional quantities. They also lead to a low-cardinality basis $\cardP$ so that the conditioning of the approximation method remains good, in the sense that $\cardP \le \Nq$, a crucial feature for deriving a good approximation from the scarce available data.

\subsection{High-Dimensional Model Representation} \label{Sec_HDMR}

An efficient alternative to these tensored-structure formats for representing high-dimensional quantities is discussed in \cite{Rabitz_Alis_99,Alis_Rabitz_01}. It consists in representing a quantity $u\left(\bxi\right)$ with a sum of lower-dimensional terms accounting for increasing levels of interaction between the constitutive variables:
\be
u\left(\bxi\right) = f_\emptyset + \sum_{i=1}^\Nd{f_i\left(\xi_i\right)} + \sum_{\substack{i,j=1, \\ j > i}}^\Nd{f_{i j}\left(\xi_i, \xi_j\right)} + \ldots + f_{1 2 \ldots \Nd}\left(\xi_1, \ldots, \xi_\Nd\right) = \sum_{\bgamma \subseteq \left\{1, \ldots, \Nd\right\}}{f_\bgamma}, \label{HDMR_full}
\ee
where $f_\bgamma$ are functions of $\randomspace$ and depend only on a subset of variables $\bxi_\bgamma = \left\{\xi_i\right\}_{i \in \bgamma}$ and $\bgamma$ is a multi-index. This decomposition is exact, unique, and does not introduce any approximation. An important property is that the modes $\left\{f_\bgamma\right\}$ are mutually orthogonal: $\left< f_\bgamma, f_{\bgamma'} \right>_{L^2\left(\Xi, \, \measxi\right)} = 0$, $\forall \, \bgamma \ne \bgamma' \subseteq \left\{1, \ldots, \Nd\right\}$.
The zero-th order term $f_\emptyset$ accounts for the mean and is invariant across the entire domain $\Xi$, while the other modes are zero-mean:
\be
f_\emptyset = \left< u \right>_{L^2\left(\Xi, \, \measxi\right)}, \qquad \left< f_\bgamma \right>_{L^2\left(\Xi, \, \measxi\right)} = 0, \quad \forall \, \bgamma \subseteq \left\{1, \ldots, \Nd\right\} \backslash \emptyset.
\ee

The rationale behind the expected success of this so-called \emph{High Dimensional Model Representation} (HDMR) is that many quantities of interest exhibit a significant dependence on low-dimensional groups of variables only, hence having negligible high order interaction decomposition terms. This leads to an efficient approximation of $u$ with only a low $\Ninter$-order HDMR: $u\left(\bxi\right) \approx  \sum_{\bgamma \subseteq \left\{1, \ldots, \Nd\right\}}{f_\bgamma\left(\bxi_\bgamma\right)}$, $|\bgamma| \le \Ninter$. We denote $\basissetf$ the set of retained modes, $\basissetf := \left\{\bgamma \subseteq \left\{1, \ldots, \Nd\right\}; |\bgamma| \le \Ninter \right\}$.

Functions $\left\{f_\bgamma\right\}$ are evaluated with the application of a set of commuting projections $\left\{\proj_i\right\}$ onto the output $u$. The projection $\proj_i$ eliminates the effect of variable $\xi_i$ while leaving the effect of the others unchanged. Letting $\proj_\emptyset$ be the identity operator on $\randomspace$, we define $\displaystyle \proj_{\boldsymbol{\eta}} = \prod_{i \in {\boldsymbol{\eta}}}{\proj_i}$, $\forall \, \boldsymbol{\eta} \subseteq \left\{1, \ldots, \Nd\right\}$. Functions $\left\{f_\bgamma\right\}$ can then be written, \cite{Kuo_etal_09},
\be
f_{\bgamma \subseteq \left\{1, \ldots, \Nd\right\} \backslash \emptyset} = \proj_{\left\{1, \ldots, \Nd\right\} \backslash \bgamma} \, u - \sum_{\bgamma' \subsetneq \bgamma}{f_{\bgamma'}} = \sum_{\bgamma' \subseteq \bgamma}{\left(-1\right)^{|\bgamma| - |\bgamma'|} \, \proj_{\left\{1, \ldots, \Nd\right\} \backslash \bgamma'} \, u}, \qquad f_\emptyset = \proj_{\left\{1, \ldots, \Nd\right\}} \, u.
\ee

Defining projections as $\proj_i \, u\left(\bxi\right) = \int_{\Xi_i}{u\left(\xi_1, \ldots, \xi_{i-1}, \zeta', \xi_{i+1}, \ldots, \xi_\Nd\right) \, \d\measbrut\left(\zeta'\right)}$, the measure $\measbrut$ determines the form of the projection. A popular choice consists in using $\measbrut = \measxioneD_i$ so that the \emph{Analysis of Variance} (ANOVA) decomposition is obtained. An example of application of the HDMR representation to the approximation of a random quantity is presented in \cite{Ma_Zabaras_10}.

\begin{remark}
These different functional representations are not totally distinct. For instance, the PC basis defined in Eq.~\eqref{PC_basis_span} can also be interpreted as a particular case of both HDMR and tensor-based expansion. For illustration, consider the following PC basis approximation space $\randomspaceNo = \span\left(\left\{\psi_{1} \left(\equiv 1\right), \psi_{2}\left(\xi_1\right), \psi_{2}\left(\xi_2\right), \psi_3\left(\xi_1\right), \psi_2\left(\xi_1\right) \, \psi_2\left(\xi_2\right), \psi_{3}\left(\xi_2\right)\right\}\right)$. This corresponds to a HDMR representation with $\Ninter = 2$ and $f_\emptyset \in \span\left(\psi_{1}\right)$, $f_1 \in \span\left(\psi_{2}\left(\xi_1\right), \psi_{3}\left(\xi_1\right)\right)$, $f_2\in \span\left(\psi_{2}\left(\xi_2\right), \psi_{3}\left(\xi_2\right)\right)$, $f_{1 2}\in \span\left(\psi_{2}\left(\xi_1\right) \psi_{2}\left(\xi_2\right) \right)$. Further, this can also be reformatted in a $\nr = 3$-rank CP format, say with $f_{1,1} \in \span\left(\psi_{1}\right)$, $f_{2,1} \in \span\left(\psi_{1}, \psi_{2}\left(\xi_2\right), \psi_{3}\left(\xi_2\right)\right)$, $f_{1,2} \in \span\left(\psi_{2}\left(\xi_1\right)\right)$, $f_{2,2} \in \span\left(\psi_{1}, \psi_{2}\left(\xi_2\right)\right)$, $f_{1,3} \in \span\left(\psi_{3}\left(\xi_1\right)\right)$ and $f_{2,3} \in \span\left(\psi_{1}\right)$.
\end{remark}

\section{Quantifying uncertainty of scattered data} \label{methodology_section}

\subsection{Setting up the stage}

In the following, we will consider that the quantity of interest $u$ is a scalar-valued random field, indexed by space and/or time $\bx \in \mathbb{R}^{\Ndx}$ and depending on a set of random variables $\bxi \in \mathbb{R}^{\Nd}$. To approximate it, the only available piece of information is a collection of scattered samples $\left\{\bxq, \bxiq, \uq\right\}_{q=1}^\Nq$. In case these data come from an experimental context, the coordinates $\bxiq$ are not directly measurable. They are then inferred from auxiliary observations and depend on the modelization.\footnote{For instance, in a fluid flow, the Reynolds number may be uncertain and modeled as a random variable parameterized by $\xi_i$. The value of $\xi_i$ in each sample $\bxiq$ is then auxiliary deduced from the measurement of the flow velocity $V$ and the model $V\left(\xi_i\right)$.}
Since the underlying random quantity $u$ is only known through these samples, no governing equation for the QoI can be exploited and, say, Galerkin projection-based weak-formulation methods cannot be employed. Further, these samples are scattered and do not follow a deterministic rule so that no deterministic sampling strategy can be assumed. Quadrature-based techniques can then not be applied either and one has to resort to regression to estimate the coefficients of the approximation in the retained basis $\left\{\psi_\alpha\right\}$. Standard $L^2$-regression solves Eq.~\eqref{LS} which is only well-posed for a matrix $\matPsi$ such that $\matPsi^T \, \matPsi$ is invertible so that it requires the number of observations to be larger than the cardinality of the approximation basis, $\Nq \ge \cardP$.

The choice of a good approximation basis in a general setting largely remains an open question. On one hand, if one is given a dictionary of approximation functions, {\apriori} selecting the best terms so that they can be evaluated from the data is a combinatorial optimization problem which algorithmic complexity quickly becomes intractable when the size of the dictionary grows. On the other hand, dictionary-learning techniques require a training while availability of an independent training set cannot be assumed here.

\begin{figure}[!ht]
\begin{center}
  \includegraphics[angle = 0, width = \textwidth, draft = false]{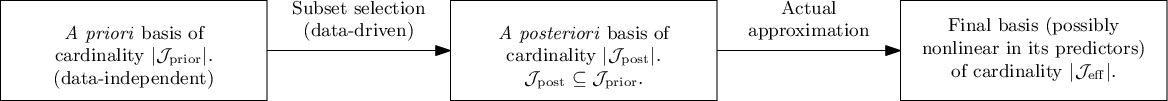}
  \caption{Sketch of the solution method.}
  \label{Fig_sketch}
\end{center}
\end{figure}

The proposed approach is as follows. We separate the determination of an efficient representation format from the evaluation of the coefficients. We first choose an {\apriori} general format for the approximation of $u$, section \ref{Aprioribasis}. The selection of particular terms to be included in the approximation basis is left to a dedicated subset selection procedure which will further refine the approximation basis and make it as tight as possible, section \ref{subset_selection_sec}. A good {\apriori} basis is motivated by results from Compressed Sensing which show that the number of samples necessary for accurately selecting the dominant basis functions of a $K$-sparse QoI (\ieLM, having $K$ non-zero coefficients in the retained approximation basis) varies as $K \, \log\left(\cardP\right)$, \cite{Candes_Romberg_06_InvPb}, illustrating the fact that it becomes increasingly difficult to select the best terms when the size $\cardP$ of the {\apriori} dictionary increases.
The subset selection hence produces an {\aposteriori} basis suitable for the data at hand. However, this basis is \emph{linear} in its predictors as required by the selection method. To circumvent this limitation, the {\aposteriori} basis is used as a skeleton only, of the best structure, and the final approximation of the QoI is evaluated with a different basis, of the same skeleton, but possibly nonlinear in its predictors, section~\ref{comp_lambda}. A sketch of the solution method is shown in Fig.~\ref{Fig_sketch}.

%
%

\subsection{A priori choice of representation of a random variable} \label{Aprioribasis}
%

We first focus on approximating a random variable and will discuss approximation of a more general random process in section~\ref{KLapprox}. The QoI is hence here a random variable $u\left(\bxi\right)$.

In this work, we want to take advantage of the low order interactions of constitutive variables for many quantities of practical interest as mentioned in section \ref{Sec_HDMR}. Previous works have shown evidence of this low interaction configuration in various situations, \cite{Rabitz_Alis_99,Alis_Rabitz_01, Ma_Zabaras_10}, and the QoI is hence chosen to be approximated under the HDMR form, Eq.~\eqref{HDMR_full}. An example is considered in \ref{Appendix_A} and demonstrates that a general HDMR format approximation with a tensor-based description of the interaction modes $\left\{f_\bgamma\right\}$ involved in the HDMR may compare favorably with a full tensor-based approximation in terms of required number $\cardPprior$ of basis functions for a given reconstruction accuracy, even for reasonably large dimensional problems. This motivates our choice of an HDMR format for the {\apriori}, data-independent, basis.

\subsection{Subset selection} \label{subset_selection_sec}

%
%
We now build upon from the {\apriori} basis and further improve it with an {\aposteriori}, data-driven, procedure.

\subsubsection{A direct approach} \label{direct_approach}

As discussed in section~\ref{Comput_coefs}, different techniques may be used to compute the coefficients of an approximation. In the case considered in this paper,
the available data are scarce while the cardinality $\cardPprior$ of the {\apriori} approximation basis may be large, in particular when the dimensionality $\Nd$ of the problem is large. It can then result in an ill-posed problem where one has to estimate $\cardPprior$ coefficients for each stochastic mode $\modesto_n$ from $\Nq \ll \cardPprior$ pieces of information. 
However, this situation often only reflects our lack of knowledge on the quantity at hand and how conservative this naive approximation method is. Indeed, high-dimensional problems are often intrinsically sparse and lower dimensional. In the present setting, it is likely that many dimensions actually hardly contribute to the approximation and that representing the dependence of the QoI along only a subset of the dimensions yields an acceptable accuracy. In our {\apriori} HDMR representation, it means that many interaction modes $\left\{f_\bgamma\right\}$ can be discarded without significantly affecting the accuracy. The challenge for an efficient solution method is then to reveal and exploit the low-dimensional manifold onto which a good approximation of the solution lies.
As an illustration, if $u\left(\bxi\right) = g\left(\xi_i\right)$ was depending only on one dimension $i$, $i \in \left\{1, \ldots, \Nd\right\}$, information theory allows to show that one only requires $m+1+\lceil\log_2 \Nd\rceil$ function evaluations to approximate a sufficiently smooth function $g \in C^s$, having $s$ continuous derivatives, so that $\left\|u - \uhat\right\|_{C\left(\Xi_i\right)} \le a \, h^s$, $h := 1/m$, where $a \ge 0$ is related to a norm of $g$, \cite{deVore_etal_11}. This number of samples actually is directly related to the number of information bits required to represent the integer $i \in [1, \Nd]$.


While determining which interaction modes are dominant is an NP-hard problem in general, recent results have shown that a good estimation of the best subset can be obtained as the solution of a convex optimization problem. In particular, the LASSO formulation, \cite{Tibshirani_96}, has been proved effective. One of its formulations, referred to as \emph{Basis Pursuit Denoising}, writes:
\be
\bcoef = \argmin_{\bcoeft \in \mathbb{R}^\cardPprior} \left\|\bcoeft\right\|_1 \quad \mathrm{s.t.} \quad \left\|\bu - \matPsi \, \bcoeft\right\|_2 \le \epsilon, \label{BPD}
\ee
with $\matPsi$ the matrix of evaluations of the approximation basis and $\epsilon$ the approximation residual. Efforts from the signal processing community, where the theory supporting these results is termed \emph{Compressed Sensing}, have demonstrated its good recovery properties in the case where $\Nq < \cardPprior$, \egLM, \cite{Chen_Donoho_Saunders_99,Candes_Tao_04, Donoho_06}. In particular, this formulation achieves provable and robust recovery bounds.\footnote{For a sufficiently incoherent set of approximation and test functions, a $K$-sparse solution $\bcoef$ to Eq.~\eqref{BPD} satisfies, \cite{Cai_etal_10}, $\displaystyle \left\|\bcoef^\star - \bcoef\right\|_2 \lessapprox h \, \left(\epsilon + \left\|\bcoef^\star - \bcoef_K^\star\right\|_1 / \sqrt{K}\right)$, where $h > 0$ is a constant depending on the set of approximation and test functions and $\bcoef_K^\star$ is the $K$-term approximation of $\bcoef^\star$ given by an oracle, \ieLM, it is the best $K$-term approximation of $\bcoef^\star$ if one 
was given full knowledge of it.}

The Compressed Sensing technique was proved very effective and is now being applied in many areas, including Uncertainty Quantification, \cite{Doostan_Owhadi_11, Mathelin_Gallivan_12}. However, standard implementations of the algorithm require the sensing matrix $\matPsi$ to be available. This bears an intrinsic limitation when it comes to high-dimensional problems as it requires the use of the whole dictionary at once from which to select the basis functions associated with the dominant coefficients.
While effective, this approach is not deemed tractable for high-dimensional problems, neither in terms of storage requirement nor CPU burden.

\subsubsection{A progressive selection}

To circumvent the issues identified above, we here use a bottom-to-top approach which achieves a forward stagewise regression by progressively revealing important basis functions. Introduced by \cite{Efron_etal_04, Hastie_etal_09}, the Least Angle Regression Selection (LARS) technique relies on analytical solutions to speed-up computations and essentially follows the piecewise linear regularization path of the LASSO.\footnote{In a nutshell, it consists in selecting, from the {\apriori} set $\basissetprior$, the predictor (approximation function) which is most correlated with the current residual, move this predictor to the active set $\basissetpost $, compute the increment solution vector by minimizing the residual $L^2$-norm and follow the descent direction along the increment vector until a predictor from the inactive set becomes as correlated with the residual as those from the active set. The whole process is then repeated and allows to sequentially build the optimal subset of approximation functions by 
exploring the Pareto front defined by the competition between the two terms of the unconstrained formulation of the optimization problem of Eq.~\eqref{BPD}.} One advantage of LARS over other techniques is that the potential dictionary is never stored nor used as a whole. A LARS approach in the UQ framework was also considered in \cite{Blatman_Sudret_2011}.

We consider the following polynomial approximation $\ft_{\bgamma}$ of $f_{\bgamma}$:
\be
f_{\bgamma}\left(\left\{\xi_i\right\}_{i \in \bgamma}\right) \approx \ft_{\bgamma}\left(\left\{\xi_i\right\}_{i \in \bgamma}\right) := \sum_{\balpha, \, |\balpha| \le \NoLARS}{\coef_{\bgamma, \balpha} \, \psi_\balpha\left(\left\{\xi_i\right\}_{i \in \bgamma}\right)}, \qquad \psi_\balpha = \prod_{i \in \bgamma}{\psi_{\alpha_i}\left(\xi_i\right)}, \label{HDMR_linearterm}
\ee
with $\balpha = \left(\alpha_i, i \in \bgamma\right)$, $\alpha_i \in \left\{1, \ldots, \Not\right\}$. Interaction modes $\left\{f_\bgamma\right\}$ are then approximated in $\mathbb{P}_{\NoLARS}$, the space of polynomials with maximum total degree $\Not$, by modes $\left\{\ft_\bgamma\right\}$ linear in their coefficients.

In the present framework, the HDMR approximation format naturally leads to \emph{groups} of predictors whose importance in describing the QoI $u$ follows a similar trend. These groups are defined by the subsets $\left\{\group_{\bgamma}\right\}$ of predictors which belong to a given interaction mode $f_{\bgamma}$, $\group_{\bgamma} = \left\{\psi_\balpha\left(\left\{\xi_i\right\}_{i \in \bgamma}\right)\right\}$, and are likely to be strongly correlated. For instance, if the QoI exhibits a strong dependence on a given dimension $\xi_j$, one then wants to incorporate the whole set of predictors $\left\{\psi_{\balpha}\left(\left\{\xi_i\right\}_{i \in \bgamma}\right)\right\}$, $\bgamma: \, j \in \bgamma$ without evaluating their relevance individually. One then looks for an approximation which is sparse at the level of groups of functions. Note that grouping predictors significantly alleviates the computational cost associated with the subset selection as further discussed in section~\ref{Complexity_sec}.

It is important to recall that this approximation format is made only for the subset selection step and is independent of the format the QoI will finally be approximated in.
The selection of groups reduces to selection of interaction modes $f_{\bgamma}$ and leaves the possibility for using different formats between the subset selection step and the coefficients evaluation step: an interaction mode found to be dominant is incorporated to the active dictionary $\basissetfpost$ independently of the way its contribution to the approximation of $u$ is actually determined in the end. Indeed, since the LARS technique only applies to predictors \emph{linear} in their coefficients, an approximation $\ft_{\bgamma}$ of the form \eqref{HDMR_linearterm} is suitable for the selection of the dominant groups. However, the final approximation $\fhat_\bgamma$ of the retained $f_\bgamma$ may rely on predictors \emph{nonlinear} in their coefficients: the subset selection step only serves to determine which interaction modes will be considered in the {\aposteriori} approximation basis, the `skeleton' $\left\{f_{\bgamma}: \, \bgamma \in \basissetfpost\right\}$.

The selection is made using a modified LARS approach and the following optimization problem is solved:
\be
\bcoef = \argmin_{\bcoeft \in \mathbb{R}^\cardP} \left\|\bu - \matPsi \, \bcoeft\right\|_2^2 + \tau \, \sum_{\bgamma \in \basissetf}{\left\|\bcoeft_{\bgamma}\right\|_{K_\bgamma}}, \label{gLARS}
\ee
with $\tau > 0$ the regularization parameter and $\left\|\cdot\right\|_{K_\bgamma}$ a norm induced by a positive definite matrix $K_\bgamma$. All predictors within a group $\bgamma$ are here weighted similarly so that we use a scaled identity matrix $K_\bgamma = \Ident_{\cardPgroup } / \cardPgroup$, $\forall \, \bgamma \in \basissetf$. The regularization term is a combination of $L^2$- and $L^1$- norms and penalizes the $L^1$-norm of the `group' vector to promote a collective behavior: either a group is basically active (non-zero $K_\bgamma$-norm) or inactive, essentially disregarding the detailed behavior within the group. This group LARS (gLARS) strategy was first proposed in \cite{Yuan_Lin_06} and the algorithm presented in \cite{Xie_Zeng_10} was modified to solve the optimization problem \eqref{gLARS}.
set of dominant modes $\left\{f_\bgamma\right\}$ is first determined by the gLARS approach with a low approximation order $\NoLARS$ and the basis is subsequently further refined by a LARS step, using $L^1$-regularization, onto these selected modes only now approximated with a higher $\NoLARS$ for improved accuracy.

%
\subsection{Functional spaces for the final approximation basis}  \label{comp_lambda}

We now discuss the general methodology for approximating a random variable $u\left(\bxi\right)$, from a finite set of its realizations. An {\apriori} choice of representation format was first made, section \ref{Aprioribasis}, and was adjusted based on the data through the subset selection procedure, the {\aposteriori} step, section \ref{subset_selection_sec}. This has selected a set of groups, or interaction modes, $\left\{f_\bgamma\right\}_{\bgamma \in \basissetfpost}$ deemed to most contribute to the HDMR representation of the QoI $u$. The actual approximation of $u$ will rely on these selected groups but does not bear restriction on the linearity w.r.t. the coefficients so that different suitable formats, possibly nonlinear, can then be considered.

Many possibilities exist to determine an approximation of $\left\{f_\bgamma, \bgamma \in \basissetfpost\right\}$ in a polynomial space, \egLM, maximum partial degree, maximum total degree, hyperbolic cross, etc. For sake of simplicity, the space $\totpolspace$ of polynomials with maximum total degree $\No$ is retained as a reasonable compromise between cardinality $\cardPgroup$ and expected accuracy of the approximation $\fhat_{\bgamma}$:
\bea
&& f_{\bgamma}\left(\left\{\xi_i\right\}_{i \in \bgamma}\right) \approx \fhat_{\bgamma}\left(\left\{\xi_i\right\}_{i \in \bgamma}\right) = \sum_{\balpha, |\balpha| \le \No}{\coef_{\bgamma, \balpha} \, \psi_\balpha\left(\left\{\xi_i\right\}_{i \in \bgamma}\right)}, \qquad \psi_\balpha\left(\left\{\xi_i\right\}_{i \in \bgamma}\right) = \prod_{i \in \bgamma}{\psi_{\alpha_i}\left(\xi_i\right)}, \nonumber \\
&& \balpha = \left(\alpha_i, i \in \bgamma\right), \qquad \alpha_i \in \left\{1, \ldots, \No\right\}, \qquad 1 \le |\bgamma| \le \NinterPC \le \min\left(\Ninter, \No\right). \label{HDMR_gal}
\eea

The cardinality associated with this approximation of $f_{\bgamma}$ at a given iteration level $l = |\bgamma|$ is $\cardPgroup = \No! / \left(l! \, \left(\No - l\right)!\right)$ and usually provides an accurate approximation with a low number of coefficients for low dimensions $|\bgamma|$.

When the dimension $|\bgamma|$ increases, the number of terms in $\fhat_{\bgamma}$ decreases and eventually degenerates for $|\bgamma| > \No$. For modes of interaction order highe than a prescribed threshold $\NinterPC$, a low-rank canonical decomposition is instead considered:
\be
f_{\bgamma}\left(\left\{\xi_i\right\}_{i \in \bgamma}\right) \approx \fhat_{\bgamma}\left(\left\{\xi_i\right\}_{i \in \bgamma}\right) = \sum_{r=1}^\nr \, \prod_{i \in \bgamma}{\sum_{\alpha=1}^\No{\coef_{\bgamma, \alpha}^{r, i} \, \psi_\alpha\left(\xi_i\right)}}, \qquad \NinterPC < |\bgamma| \le \Ninter \le \Nd. \label{HDMR_term}
\ee

The maximum number of modes at a given interaction level $l$ is $\Nd! / \left(\left(\Nd-l\right)! \, l!\right)$. Relying on an approximation in $\totpolspace$ for interaction modes of order $|\bgamma| \le \NinterPC$ and on low-rank approximation for higher interaction order modes, with maximum rank $\nr$, the total cardinality of this approximation format is bounded from above by
\be
\cardPcur \leq \sum_{l=0}^{\NinterPC}{\frac{\Nd! \, \No!}{\left(\Nd - l\right)! \, \left(l!\right)^2  \, \left(\No - l\right)!}} + \sum_{l=\NinterPC+1}^\Ninter{\frac{\Nd!}{\left(\Nd - l\right)! \, l!} \, \nr \, l \, \No}.
\ee

\subsection{Algorithm for approximating a random variable} \label{Algo_lambda}

We will denote $\basissetfcur$ the set of modes $\left\{\fhat_{\bgamma}\right\}_{\bgamma \subseteq \basissetfpost}$ finally considered for the approximation of $u$ and $\basissetcur$ the set of associated predictors $\left\{\psi_\balpha\right\}$. The interaction modes are estimated sequentially. Once a new mode is evaluated, the whole approximation may be updated by reevaluating the coefficients of the predictors $\left\{\psi_\balpha\right\}$ already evaluated of the current evaluation set $\basissetfcur \in \basissetfpost$. Let $\bres = \left(\res^{(1)} \, \ldots \res^{(\Nq)}\right)^T$ be the residual vector after basis functions $\fhat_{\bgamma}, \bgamma \in \basissetfcur$ have been evaluated. The coefficients involved in the next mode $\fhat_{\bgamma}, \bgamma \in \basissetfpost \backslash \basissetfcur$ to be evaluated are then determined. If $\bgamma$ is such that $|\bgamma| \le \NinterPC$, they are computed from the following system of equations\footnote{While not found necessary here, the solution of 
the least squares problem may be regularized by adding a generic term of the form $\beta \: \left\| L \, \bcoeft\right\|_2$. A typical choice is $L = I_{|\basisset_{\bgamma}|}$ but one may also want to consider non-diagonal matrices $L$.}:
\be
\left\{ \,
\bcoef_{\bgamma, \cdot} = \argmin_{\bcoeft \in \mathbb{R}^{|\basisset_{\bgamma}|}}{\left\|\bres - \Psi \, \bcoeft\right\|_2}, \qquad  \forall \, i \in \bgamma, \quad \bgamma \subseteq \basissetfpost, \quad |\bgamma| \le \NinterPC, \right. \label{lambda_system_linear}
\ee
with $\bcoef_{\bgamma, \cdot} = \left(\bcoef_{\bgamma, \alpha_i}, i \in \bgamma\right)^T$ and
\bea
\res^{(q)} & = & \uq - \sum_{\bgamma' \subseteq \basissetfcur \backslash \bgamma}{\fhat_{\bgamma'}\left(\left\{\xiq_i\right\}_{i \in \bgamma'}\right)}, \qquad \bres = \left(\res^{(1)} \, \ldots \, \res^{(\Nq)}\right)^T, \nonumber \\
\Psi_{q \balpha} & = & \psi_\balpha\left(\left\{\xiq_i\right\}_{i \in \bgamma}\right), \qquad \Psi = \left[\Psi_{q \balpha}\right].
\eea

To solve for the coefficients associated with predictors nonlinear in their coefficients, an Alternate Least Squares (ALS) approach is used, reformulating the nonlinear problem into a set of coupled linear equations:
\be
\left\{ \,
\bcoef_{\bgamma, \cdot}^{r, i} = \argmin_{\bcoeft \in \mathbb{R}^\No}{\left\|\bres_i - \Psi \, \bcoeft\right\|_2}, \qquad \forall \, i \in \bgamma \subseteq \basissetfpost, \quad \NinterPC < |\bgamma| \le \Ninter, \right. \label{lambda_system}
\ee
with $\bcoef_{\bgamma, \cdot}^{r, i} = \left(\bcoef_{\bgamma, 1}^{r, i} \, \ldots \, \bcoef_{\bgamma, \No}^{r, i} \right)^T$ and
\bea
\res_i^{(q)} & = & \uq - \sum_{\bgamma' \subseteq \basissetfcur}{\fhat_{\bgamma'}\left(\left\{\xiq_i\right\}_{i \in \bgamma'}\right)} - \sum_{r'=1}^{r-1}{\prod_{i' \in \bgamma}{\sum_{\alpha'=1}^{\No}{ \coef_{\bgamma, \alpha'}^{r', i'} \, \psi_{\alpha'}\left(\xiq_{i'}\right) }}}, \nonumber \\
\Psi_{q \alpha} & = & \psi_{\alpha}\left(\xiq_{i}\right) \, \prod_{i' \in \bgamma, i' \ne i}{\sum_{\alpha'=1}^{\No}{ \coef_{\bgamma, \alpha'}^{r, i'} \, \psi_{\alpha'}\left(\xiq_{i'}\right)}}, \qquad \Psi = \left[\Psi_{q \balpha}\right].
\eea

This whole step is embedded in a loop over the modes $f_{\bgamma}, \bgamma \in \basissetfpost$ retained by the subset selection procedure. The cross-validation error (CV$\varepsilon$) is estimated from $\Nqhat$ validation samples $\left\{\bxiqh, \uqh\right\}_{\hat{q}=1}^{\Nqhat}$ independent from the $\Nq$ samples of the training set.\footnote{A ratio $\Nqhat / \Nq \simeq 1 / 2$ is typically accepted as a reasonable splitting of the set of samples. We here use the simplest cross-validation method but more sophisticated techniques ($k$-fold, Leave-One-Out, etc.) are available, see for instance \cite{Hastie_etal_09}. While more accurate, they are significantly more computationally expensive.} If the cross-validation error  has increased over the last two loops, the approximation basis is likely to have become too large w.r.t. the available data and iterations are stopped. The retained basis is then the one that has led to the lowest CV$\varepsilon$. On the other hand, if CV$\varepsilon$ keeps decreasing, the 
next interaction mode as selected by the subset selection step is considered and added to the current active set $\basissetfcur$ and the whole iteration is carried-out.
Once the approximation is determined, the coefficients are updated with the same sequential technique using both the training and the validation points, $\Nq + \Nqhat$. The approximation accuracy is estimated by the relative $L^2$-norm $\varepsilon$ of the approximation error estimation evaluated from a $\Nqt$-point test set $\left\{\left(\xqt, \bxiqt, \uqt\right)\right\}_{\tilde{q}=1}^{\Nqt}$, independent from the training set:
\be
\varepsilon^2 := \left\| \bu - \buhat \right\|^2_2 / \left\| \bu \right\|^2_2, \qquad \bu = \left(u^{(1)} \, \ldots \, u^{(\Nqt)}\right), \qquad \buhat = \left(\uhat^{(1)} \, \ldots \, \uhat^{(\Nqt)}\right). \label{sto_error_def}
\ee
The global methodology is summarized in Algorithm~\ref{Algo_sumup}. Statistical moments can be readily evaluated from the present HDMR of the QoI, see \ref{Sobol_sec}.

\begin{algorithm}
\caption{Sketch of the solution method for approximating a random variable $u\left(\bxi\right)$}
\label{Algo_sumup}
\begin{algorithmic}[1]
  \STATE Select an {\apriori} basis in HDMR format. Choose $\No$, $\Ninter$, $\NinterPC$, $\nr$ and $\NoLARS$. Initialize $\bres = \left(u^{(1)} \ldots u^{(\Nq)}\right)^T$.
  \STATE \textbf{Subset selection step}. Solve the LASSO optimization problem with the gLARS algorithm $\lra$ sequence of {\aposteriori} approximation bases indexed by $s$ with ordered groups $\basissetfpost = \left\{\bgamma^{(s)}\right\}$. Initialize $s$ and $\basissetfcur$: $s \leftarrow 0$, $\basissetfcur \leftarrow \emptyset$. \label{subsetselstep}
  \STATE \textbf{Solve the approximation problem}: \label{solvestep}
      \REPEAT
	  \STATE $s \leftarrow s+1$.
          \STATE Consider the next mode $\fhat_{\bgamma^{(s)}}$ from the set $\basissetfpost$ selected in (\ref{subsetselstep}): $\basissetfcur \leftarrow \basissetfcur \bigcup \bgamma^{(s)}$.
	  \STATE Solve for the approximation coefficients $\left\{\bcoef_{\bgamma}\right\}_{\bgamma \in \basissetfcur}$ by alternately solving for the coefficients of modes $\left\{\fhat_{\bgamma}\right\}_{\bgamma \in \basissetfcur}$, Eqs.~(\ref{lambda_system_linear}, \ref{lambda_system}). [\textbf{Update step}] \label{Updateornot}
	  \STATE Estimate the cross-validation error CV$\varepsilon$ and evaluate the current approximation $\buhat = \left(\uhat^{(1)} \ldots \uhat^{(\Nq)}\right)^T$.
	  \STATE Update the residual $\bres \leftarrow \bu - \buhat$.
      \UNTIL {CV$\varepsilon$ has increased over the last two passes $s$ and $s-1$.}
  \STATE $\basissetfcur \leftarrow \basissetfcur \backslash \left\{\bgamma^{(s)}, \bgamma^{(s-1)}\right\}$.
  \STATE Update the coefficients $\left\{\bcoef_{\bgamma}\right\}_{\bgamma \in \basissetfcur}$ of the retained modes with the extended set of data $\left\{\bxiq, \uq\right\}_{q=1}^{\Nq + \Nqhat}$. It finally yields $\uhat\left(\bxi\right)$ expressed in the basis $\left\{\fhat_{\bgamma}\right\}_{\bgamma \in \basissetfcur}$. \label{updatestep}
\end{algorithmic}
\end{algorithm}

\subsection{Robust estimation} \label{Noise_section}

An important concern when deriving a methodology is the robustness w.r.t. noise and a more robust alternative to the methodology discussed so far is now presented.

To evaluate the approximation coefficients once an approximation basis is determined from the subset selection step, a standard approach is to minimize a norm between target observations and reconstructed approximation as done in the previous section, Eqs. \eqref{lambda_system_linear} and \eqref{lambda_system}: the approximation coefficients of a given mode $\fhat_{\bgamma}$ are basically given by $\bcoef_{\bgamma} = \argmin_{\bcoeft \in \mathbb{R}^{\cardPgroup}} \left\|\bres - \matPsi \, \bcoeft\right\|_2$, with $\matPsi \in\mathbb{R}^{\Nq \times \cardPgroup}$ the matrix of the $\bgamma$-group predictors evaluated in $\left\{\bxiq\right\}$ and $\bres$ the target residual vector. The solution to this least squares problem is equivalently obtained from
\be
\left\{\bcoef, \Delta \bres\right\} = \argmin_{\bcoeft \in \mathbb{R}^{\cardPgroup}} \left\|\widetilde{\Delta \bres}\right\|_F \quad \mathrm{s.t.} \quad \bres + \widetilde{\Delta \bres} = \matPsi \, \bcoeft, \label{wTLS_1}
\ee
which minimizes the Frobenius norm of the residual vector. This implicitly assumes no error in the coordinates $\left\{\bxiq\right\}$ at which the target is evaluated. For instance, these coordinates may be known as the solution of auxiliary inference problems. This brings errors so that the actual coordinates vector is only estimated with an error $\Delta \bxiq$. Since $\matPsi$ depends on $\bxi$, an error predictor matrix $\Delta \matPsi\left(\bxi, \Delta \bxi\right) := \matPsi\left(\bxi + \Delta \bxi\right) - \matPsi\left(\bxi\right)$ arises and the estimation problem \eqref{wTLS_1} then rewrites as a Total Least Squares problem, \cite{Golub_vanLoan_12}:
\be
\left\{\bcoef, \Delta \matPsi, \Delta \bres\right\} = \argmin_{\bcoeft \in \mathbb{R}^{\cardPgroup}} \left\|\widetilde{\Delta \matPsi} \, \widetilde{\Delta \bres}\right\|_F \quad \mathrm{s.t.} \quad \bres + \widetilde{\Delta \bres} = \left(\matPsi + \widetilde{\Delta \matPsi}\right) \, \bcoeft. \label{wtls_gal}
\ee

The realizations of the error in the data $\left\{\Delta \bxiq\left(\theta\right), \Delta \resq\left(\theta\right)\right\}$ are modeled to follow the distribution of zero-mean \textit{iid} variables. Further, predictors may be correlated:
\be
\esp_\theta\left[\left(\Delta \psi_\balpha - \esp_\theta\left[\Delta \psi_\balpha\right]\right) \, \left(\Delta \psi_{\balpha'} - \esp_\theta\left[\Delta \psi_{\balpha'}\right]\right)\right] \ne 0,
\ee
with $\Delta \psi_\balpha = \Delta \psi_\balpha\left(\bxiq, \Delta \bxiq\right)$. A general approach to solve the weighted Total Least Squares (wTLS) problem of Eq.~\eqref{wtls_gal} consists in the minimization of the usual weighted residual sum of squares $\rho^2$, \cite{Markovsky_vanHuffel_07}:
\be
\rho^2 := \vectorize\left(\Delta X\right)^T \, \Lambda^{-1} \, \vectorize\left(\Delta X\right), \qquad \Delta X := \left(\Delta \matPsi \, \Delta \bres\right)^T, \qquad X := \left(\matPsi \, \bres\right)^T,
\label{MLPCA_correl}
\ee
where the `vec' operator unfolds a generic $m \times n$ matrix into a $m n$ vector and $X$ is the data matrix. The covariance matrix for $X$, $\Lambda := \left<\mathrm{vec}\left(X - \left<X\right>_\Nq\right) \, \mathrm{vec}\left(X - \left<X\right>_\Nq\right)^T\right>_\Nq$, is evaluated and the minimization problem \eqref{wtls_gal} is solved using the ALS-based algorithm proposed in \cite{Wentzell_etal_97}.
%

As will be shown in the numerical experiments examples, section \ref{Robust_results}, the present total least squares formulation allows to improve the approximation quality from noisy data.

\begin{remark}
When a large amount $\Nq$ of experimental information is available, the data matrix $X \in \mathbb{R}^{\left(\cardPgroup + 1\right) \times \Nq}$ can be large. The resulting correlation matrix $\Lambda$ then has potentially very large dimensions. However, since the noise is assumed independent from one sample to another, $\Lambda$ has a block diagonal structure. Further, it is a symmetric definite positive matrix, allowing for additional reduction of the storage requirement. The structure of $\Lambda$ is then exploited in solving the weighted total least squares problem above through sparse storage and operations.
\end{remark}

\subsection{Asymptotic numerical complexity}	\label{Complexity_sec}

While the primary motivation for this work is to determine an accurate representation of a random quantity from a small set of its realizations, it is desirable that the solution method remains computationally tractable. As seen above, the algorithm for approximating a random variable is essentially two fold.

The selection process essentially consists in sequentially building a subset, section~\ref{subset_selection_sec}. Each step of the sequence involves solving a least squares problem of growing size and finding the basis function, or group of functions, within the {\apriori} set $\basissetprior$ most correlated with the current residual. The matrix of the least squares problem is $\matPsi \in \mathbb{R}^{\Nq \times \cardPpost}$, with $\cardPpost$ the cardinality of the current set of selected basis functions. The least squares problem is solved \textit{via} a QR decomposition of $\matPsi$ in $\mathcal{O}\left(\Nq \, \cardPpost^2\right)$ operations. The iterative selection process is carried-out with a growing active set $\basissetpost$ until the problem becomes ill-posed, \ieLM, until $\cardPpost$ is about $\Nq$. We use grouped LARS and denote $\cardPgroupmean$ the average cardinality of the retained group predictors, \ieLM, the average number of basis functions in the group added to the active set. The subset 
selection process retains $\nJ$ groups of variables so that the total cost associated with the least squares step of the subset selection is
\be
\mathscr{J}_{\rm LS} = \mathcal{O}\left(\sum_{s=1}^{\nJ}{\Nq \, \left(\cardPgroupmean \, s\right)^2}\right).
\ee

As groups of predictors are moved to the active set, the size of the remaining {\apriori} set decreases, $\cardPprior^{(\rm current)} \simeq \cardPprior - s \, \cardPgroupmean$. The cost associated with the evaluation of the correlation for each predictor in the inactive set is then:
\be
\mathscr{J}_{\rm correl} = \mathcal{O}\left(\sum_{s=1}^{\nJ}{\Nq \, \left(\cardPprior - s \, \cardPgroupmean\right)}\right) \propto \Nq. \label{J_correl}
\ee


In practice, the cost associated with the evaluation of the correlation of the predictors in the inactive set with the current residual dominates so that the whole cost of the subset selection finally approximates as
\be
\mathscr{J}_{\rm subsel} = \mathscr{J}_{\rm LS} + \mathscr{J}_{\rm correl} \simeq \mathcal{O}\left(\Nq \, \cardPprior \, \nJ - \Nq \, \cardPgroupmean \, \frac{\nJ \, \left(\nJ+1\right)}{2}\right). \label{theo_cost_subset}
\ee

The second step of the solution method deals with the evaluation of the approximation coefficients, sections \ref{comp_lambda}-\ref{Algo_lambda}.
The cost associated with evaluating the coefficients of a $l$-th interaction order mode, $1 \le l \le \NinterPC$, encompasses the matrix $\Psi$ assembly cost $\mathcal{O}\left(\Nq \, \No! / \left(l! \, \left(\No - l\right)!\right)\right)$ and the least squares solution $\mathcal{O}\left(\Nq \, \left(\No! / \left(l! \, \left(\No - l\right)!\right)\right)^2\right)$. Since modes $\left\{\fhat_{\bgamma}\right\}_{\bgamma \in \basissetfcur}$ already evaluated may be updated once an additional one from the selected set is considered, the total cost is the sum of an arithmetic sequence. Its exact formulation depends on the selected set and is difficult to derive in closed-form. As a simple example, updating all coefficients for each new mode $\fhat_{\bgamma}$ considered, neglecting the cost associated with first-order interaction modes and assuming only second-order interaction modes are retained in the {\aposteriori} set, an upper bound for the cost writes
\be
\mathscr{J}_\mathrm{coef} \le \sum_{s=1}^{\cardPfcur}{\left[\mathcal{O}\left(s \, \Nq \, \No^l\right) + \mathcal{O}\left(s \, \Nq \, \No^{2 l}\right)\right]}, \qquad {\rm with~} l = 2,
\ee
where $\cardPfcur$ is the number of groups finally retained for the approximation by the CV test, see Algorithm \ref{Algo_sumup}. A quantitative discussion of the numerical cost is given in section \ref{Num_cost_example} with an illustrative example.



\subsection{Approximation of a random process by a separated representation} \label{KLapprox}

The approximation of a random process, say, a space-dependent uncertain quantity $u\left(\bx,\bxi\right)$ is now considered in the form of separation of variables:
\be
u\left(\bx,\bxi\right) \approx \modex_0\left(\bx\right) + \sum_{n=1}^\KLrank{\modex_n\left(\bx\right) \, \modesto_n\left(\bxi\right)} \equiv \sum_{n=0}^\KLrank{\modex_n\left(\bx\right) \, \modesto_n\left(\bxi\right)}, \qquad \modesto_0 \equiv 1.\label{KL_approx_format}
\ee

The `spatial' modes are associated with all physical dimensions the random process may be indexed upon (space, time, \ldots) so that $\bx = \left(x_1 \, x_2 \ldots t \ldots\right) \subseteq \mathbb{R}^\Ndx$. They are defined as: $\modex_n\left(\bx\right) = \sum_{l=1}^\cardx{\coef^{(\modex)}_{l, n} \, \phi_l\left(\bx\right)}$ with $\left\{\phi_l\right\}$ a chosen truncated basis of cardinality $\cardx$. The functional form of `stochastic' modes $\left\{\modesto_n\right\}$ and their evaluation was discussed in sections \ref{Aprioribasis}-\ref{Algo_lambda}.

The spatial and stochastic modes of the approximation \eqref{KL_approx_format} are sequentially determined in turn.
Let $\left\|v\right\|_{\Nq} := \left<v, v\right>_\Nq$ be the norm induced by the data-driven inner product: $\left<\cdot, \cdot\right>_\Nq: \, \mathbb{R} \times \mathbb{R} \rightarrow \mathbb{R}, \left(v, w\right) \mapsto \left<v, w\right>_\Nq := \sum_{q = 1}^\Nq{v^{(q)} \, w^{(q)}}$. Assuming $\left\{\modesto_n\right\}$ known and projecting Eq.~\eqref{KL_approx_format} onto the space spanned by $\left\{\phi_l\right\}$, the coefficients $\left\{\coef^{(\modex)}_{l, n}\right\}_l$ of the deterministic mode $\modex_n$ are the solution of the following problem:
\bea
&& \left<u, \phi_k \, \modesto_{n}\right>_{\Nq} = \left<\sum_{n'=0}^{n-1}{\modex_{n'} \, \modesto_{n'}} + \sum_{l=1}^\cardx{\coef^{(\modex)}_{l, n} \, \phi_l} \, \modesto_{n}, \phi_k \, \modesto_{n}\right>_{\Nq}, \qquad \forall \, 1 \le k \le \cardx, \nonumber \\
\Leftrightarrow && \bcoef^{(\modex)}_{\cdot, n} = \argmin_{\bcoeft \in \mathbb{R}^{\cardx}} \left\|\bu - \sum_{n'=0}^{n-1}{\bmodex_{n'} \odot \bmodesto_{n'}} - \left(\Phi \, \bcoeft\right) \odot \bmodesto_{n}\right\|_2, \label{x_modes}
\eea
where $\Phi \in \mathbb{R}^{\Nq \times \cardx}$, $\Phi_{q l} = \phi_l\left(\bxq\right)$, $\bmodex_n = \left(\modex_n\left(x^{(1)}\right) \, \ldots \, \modex_n\left(x^{(\Nq)}\right)\right)^T$, $\bmodesto_n = \left(\modesto_n\left(\bxi^{(1)}\right) \, \ldots \, \modesto_n\left(\bxi^{(\Nq)}\right)\right)^T$ and $\odot$ is the Hadamard product.
Similarly, the stochastic mode $\modesto_n$ is evaluated by determining the set of coefficients $\left\{c^{(\modesto)}_{n}\right\}$ minimizing $\left\|\bu - \sum_{n'=0}^{n-1}{\bmodex_{n'} \odot \bmodesto_{n'}} - \bmodex_n \odot \bmodesto_n\left(\left\{c^{(\modesto)}_{n}\right\}\right)\right\|_2$ using Algorithm \ref{Algo_sumup} presented in section~\ref{Algo_lambda}.
The spatial mode $\modex_n$ is then evaluated from Eq.~\eqref{x_modes} given all the other information and the whole iteration is repeated until convergence of the pair $\left\{\modex_n, \modesto_n\right\}$. The next pair can then be determined with the same methodology with $n \leftarrow n + 1$. The algorithm is summarized in Algorithm \ref{Algo_separated}.

\begin{algorithm}
\caption{Skeleton of the solution method for approximating a random process}
\label{Algo_separated}
\begin{algorithmic}[1]
  \STATE Choose $\cardx$. Initialize $\bres = \bu$ and $\bmodesto_0 = \boldsymbol{1}$ and set $n \leftarrow 0$.
  \STATE \textbf{Solve for coefficients} $\left\{\coef^{(\modex)}_{l, n}\right\}_{l = 1}^\cardx$ \textbf{of the deterministic mode} using Eq.~\eqref{x_modes} and normalize $\modex_n\left(\bx\right) = \sum_l{\coef^{(\modex)}_{l, n} \, \phi_l\left(\bx\right)}$. \label{coef_solve}
  \STATE \textbf{Solve for the coefficients of the stochastic mode} $\modesto_n\left(\bxi\right)$ using Algorithm \ref{Algo_sumup} given $\bmodex_n$ and $\bres$. If $n = 0$, $\bmodesto_0 \leftarrow \boldsymbol{1}$.
  \item If $\left\|\modesto_n\right\|_{\Nq}$ converges, set $\bres \leftarrow \bres - \bmodex_n \odot \bmodesto_n$, and $n \leftarrow n + 1$. Otherwise, iterate in (\ref{coef_solve}).
  \STATE Iterate to step (\ref{coef_solve}) unless a termination criterion is met (for instance, $\left\|\modesto_n\right\|_\Nq$ below a threshold or maximum rank $\KLrank$ reached).
\end{algorithmic}
\end{algorithm}

\begin{remark}
If the separated approximation grows beyond a few modes, it is beneficial to update the coefficients of, say, the spatial modes for improved accuracy: solve for $\left\{\modex_0, \ldots, \modex_n\right\}$ given $\left\{\bu, \modesto_0, \ldots, \modesto_n\right\}$.
\end{remark}

\section{Numerical experiments} \label{Sec_results}

The methodology developed in the previous sections is now demonstrated on a set of examples. Different aspects of the global solution method are illustrated on a 1-D stochastic diffusion equation. A more computationally involved example is next considered with a Shallow Water problem with multiple sources of uncertainty.

\subsection{Stochastic diffusion equation} \label{StoDifEq}

We consider a steady-state stochastic diffusion equation on $\Omega \times \Xi$, $\Omega = \left[x_-, x_+\right] \subset \mathbb{R}$, with deterministic Dirichlet boundary conditions:
\be
\nabla_x \, \left(\nu\left(x, \bxi'\right) \, \nabla_x u\left(x,\bxi\right)\right) = F\left(x, \bxi''\right), \qquad u\left(x_-,\bxi\right) = u_-, \quad u\left(x_+,\bxi\right) = u_+. \label{diff_eq}
\ee

The right-hand side $F$ is a random source field and $\nu$ is a space-dependent random diffusion coefficient modeled as:
\bea
\nu\left(x, \bxi'\right) = \nu_0\left(x\right) + \sum_{k=1}^{\Nd_\nu}{\sqrt{\sigma_{\nu,k}} \, \omega_{\nu,k}\left(x\right) \, \xi_k'}, \qquad \bxi' = \left(\xi_1' \, \ldots \, \xi_{\Nd_\nu}'\right), \nonumber \\
F\left(x, \bxi''\right) = F_0\left(x\right) + \sum_{k=1}^{\Nd_F}{\sqrt{\sigma_{F,k}} \, \omega_{F,k}\left(x\right) \, \xi_k''}, \qquad \bxi'' = \left(\xi_1'' \, \ldots \, \xi_{\Nd_F}''\right), \label{nu_F_definition}
\eea
with $\nu_0\left(x\right) = 1$ and $F_0\left(x\right) = -1$ the respective mean values. The random variables $\left\{\xi_1', \, \ldots, \, \xi_{\Nd_\nu}', \, \xi_1'', \, \ldots, \, \xi_{\Nd_F}''\right\}$ are chosen mutually independent and uniformly distributed on $\left[0, 1\right]$. The spatial modes $\omega_{\nu,k}\left(x\right)$ and $\omega_{F,k}\left(x\right)$, and their associated amplitude $\sqrt{\sigma_{\nu,k}}$ and $\sqrt{\sigma_{F,k}}$, are the first dominant eigenfunctions of the following eigenproblems:
\bea
\int_\Omega{K_\nu\left(x,x'\right) \, \omega_{\nu, k}\left(x'\right) \, \d x'} = \sigma_{\nu, k} \, \omega_{\nu, k}\left(x\right), \qquad K_\nu\left(x,x'\right) = \sigma_\nu^2 \, e^{- \frac{\left(x-x'\right)^2}{2 \, L_{c, \nu}^2}},\nonumber \\
\int_\Omega{K_F\left(x,x'\right) \, \omega_{F, k}\left(x'\right) \, \d x'} = \sigma_{F, k} \, \omega_{F, k}\left(x\right), \qquad K_F\left(x,x'\right) = \sigma_F^2 \, e^{- \frac{\left(x-x'\right)^2}{2 \, L_{c, F}^2}}, \label{KL_eigen}
\eea
with $K_\nu$ and $K_F$ the correlation kernels. The random fields properties are chosen as $\sigma_\nu = 0.7$, $\sigma_F = 0.7$, $L_{c, \nu} = 0.3$, $L_{F, \nu} = 0.3$. Note that $F\left(\cdot, \bxi''\right) < 0$ a.e. and $\nu\left(\cdot, \bxi'\right) > 0$ a.e. so that the problem remains coercive.
The spectra of the operators associated with these eigenproblems are here the same and decay quickly thanks to the high correlation length as can be appreciated from Table~\ref{tab_lambda} where the dominant eigenvalues are given. The resulting problem is then anisotropic in $\Xi$ in the sense that the degree of dependence of the input random parameters along the dimensions $\left\{\xi_1, \ldots, \xi_8\right\}$ strongly varies.

\begin{table}[!ht]
\centering
\begin{tabular}{cccccccc}
\hline \hline
$\sigma_1$ & $\sigma_2$ & $\sigma_3$ & $\sigma_4$ & $\sigma_5$ & $\sigma_6$ & $\sigma_7$ & $\sigma_8$ \\
\hline
0.1815 & 0.1396 & 0.0906 & 0.0450 & 0.0236 & 0.0097 & 0.0035 & 0.0011 \\
\hline
\end{tabular}
\caption{Upper part of the spectrum of both eigenproblems~\eqref{KL_eigen}.}
\label{tab_lambda}
\end{table}


Denoting $\bxi = \left(\bxi' \, \bxi''\right) \in \mathbb{R}^\Nd$, $\Nd = \Ndnu + \NdF$, the solution $u$ is approximated in a rank-$\KLrank$ separated form: $u\left(x,\bxi\right) \approx \hat{u}\left(x,\bxi\right) = \modex_0\left(x\right) + \sum_{n=1}^{\KLrank}{\modex_n\left(x\right) \, \modesto_n\left(\bxi\right)}$. The stochastic approximation basis relies on a HDMR format with a maximum interaction order $\Ninter = 3$ and 1-D Legendre polynomials $\left\{\psi_\alpha\right\}_{\alpha=1}^{\No}$ of maximum degree $\No = 8$.

In this section, the focus is on approximating a purely random quantity, \ieLM, disregarding its spatial dependence. We then rely on samples of the solution $u\left(x, \bxi\right)$ taken at a given spatial location $x^\star$: $\left\{\uq := u\left(x^\star, \bxiq\right)\right\}_{q=1}^\Nq$.

\subsubsection{Influence of the number of samples}

\begin{figure}[!ht]
\begin{center}
  \includegraphics[angle = -90, width = 0.5\textwidth, draft = false]{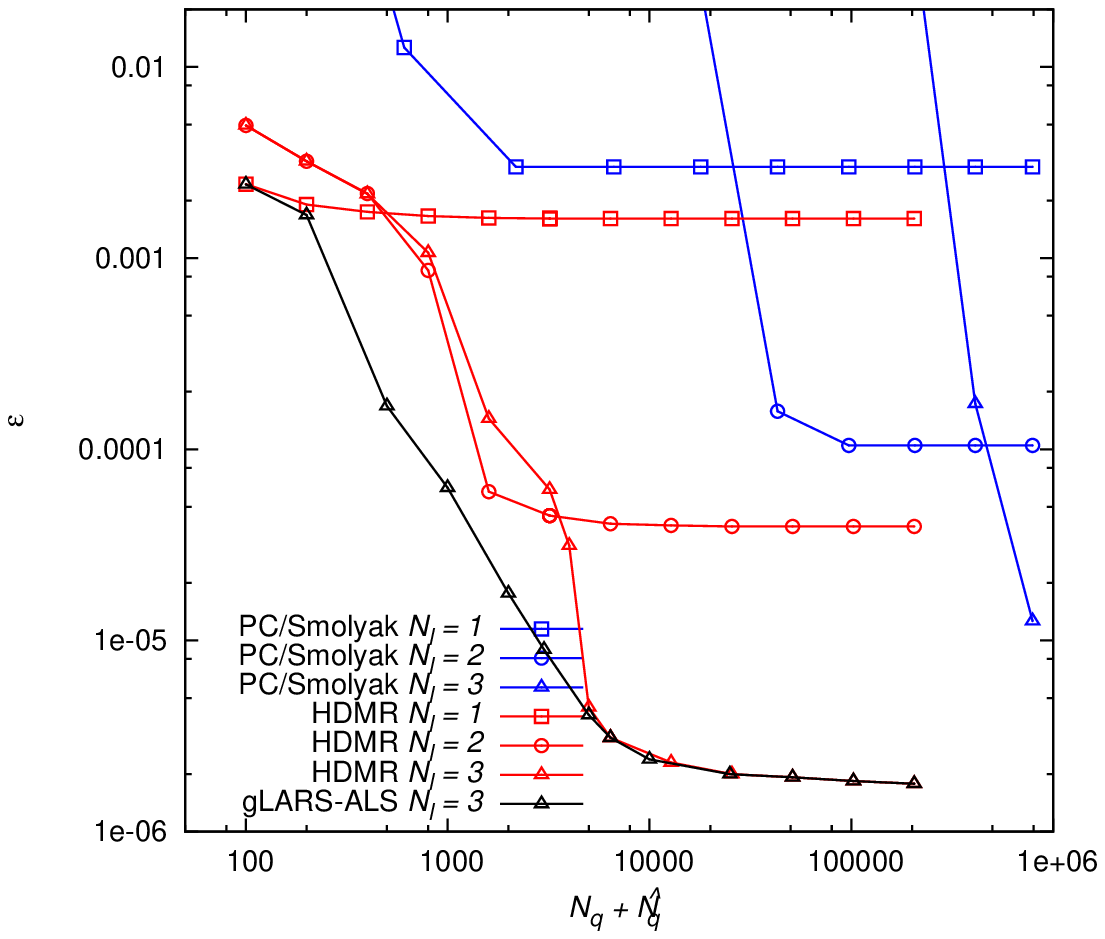}
  \caption{Convergence of the approximation with the number of samples $\Nq + \Nqhat$. Different approximation methods are compared: plain HDMR, PC/Smolyak scheme sparse grid spectral decomposition and the present gLARS-ALS. The convergence is plotted in terms of $\varepsilon$. $\Nd = 8$, $\No = 8$, $\Ninter = 3$, $\NinterPC = 3$.}
  \label{L2_MM_schemes_Nd8}
\end{center}
\end{figure}

We first focus on the achieved accuracy in the approximation with a given budget $\Nq + \Nqhat$ samples. The number of test points $\Nqt$ to estimate the approximation error $\varepsilon$, Eq.~\eqref{sto_error_def}, is chosen sufficiently large so that $\varepsilon$ is well estimated, $\Nqt = 10,000$. In Fig.~\ref{L2_MM_schemes_Nd8}, the performance of the present gLARS-ALS methodology is compared with both a plain HDMR approximation, \ieLM, with no subset selection hence considering the whole {\apriori} approximation basis, and a PC approximation with a sparse grid technique. The Smolyak scheme associated with a Gauss-Patterson quadrature rule is used as the sparse grid, with varying number of points in the 1-D quadrature rule and varying levels. The dimensionality of the stochastic space is $\Nd = 8$.

The sparse grid is seen to require a large number of samples to reach a given approximation accuracy.\footnote{Note that the plain Smolyak scheme is used here, which does not exploit anisotropy in the response surface. More sophisticated Smolyak-based approximations have been developed, see \cite{Nobile_al_07}, and are expected to provide better results.
} The HDMR-format approximation, with various interaction orders $\Ninter$, provides a better performance than PC/Smolyak but still requires more points to reach a given accuracy than the present gLARS-ALS method which performs significantly better in approximating the QoI from a given dataset. The gLARS-ALS approximation error is also seen to be smooth and monotonic when the amount of information varies. When $\Nq$ is large enough, the subset selection step becomes useless as all $\cardPprior$ terms of the {\apriori} basis can be evaluated from the large amount of information and the gLARS-ALS performance is then be similar as that of the HDMR. Note that the benefit of a subset selection step in terms of accuracy improvement increases with the dimension $\Nd$ as the size $\cardPprior$ of the potential dictionary then grows.

\subsubsection{Influence of the stochastic dimension}

\begin{figure}[!ht]
\begin{center}
  \includegraphics[angle = -90, width = 0.5\textwidth, draft = false]{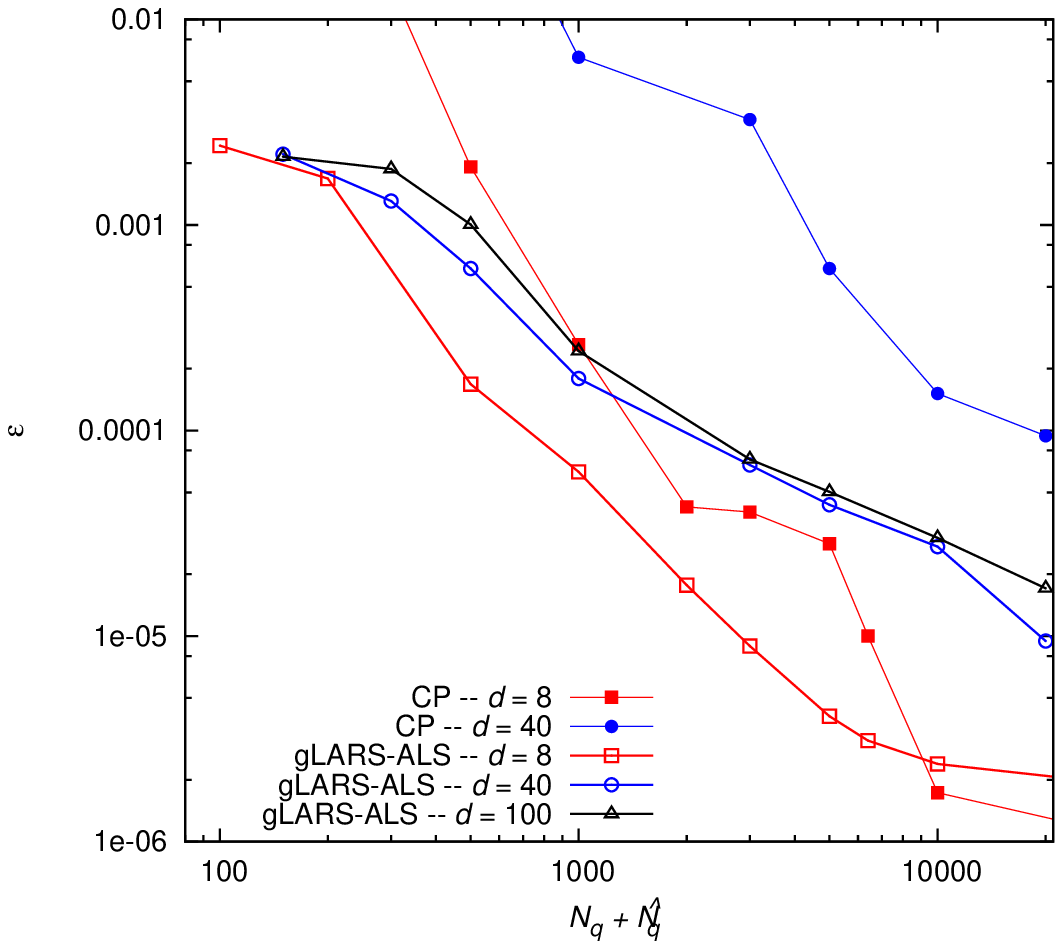}
  \caption{Convergence of the approximation with the number of samples $\Nq + \Nqhat$ and for different dimensionality of the QoI. The present gLARS-ALS approach is compared with a CANDECOMP-PARAFAC-type technique (labeled `CP').}
  \label{L2_MM_schemes_Nds}
\end{center}
\end{figure}

The approximation accuracy of the present method is now studied when the dimension of the stochastic space varies. The same problem as above is considered but with various truncation orders of the source $F$ and the diffusion coefficient $\nu$ definitions, see Eqs.~\eqref{nu_F_definition}. The solution of the diffusion problem \eqref{diff_eq} is of dimension $\Nd = \Nd_F + \Nd_\nu$ and the dimensions $\Nd_F$ and $\Nd_\nu$ are varied together, $\Nd_F = \Nd_\nu$. The resulting approximation error is plotted in Fig.~\ref{L2_MM_schemes_Nds} for different $\Nd$ when the number of available samples varies. From $\Nd = 8$ to $\Nd = 40$, the required number of points for a given accuracy is seen to increase significantly, between a 2- and a 10-fold factor. However this is much milder than the increase in the potential approximation basis cardinality, \ieLM, if not subset selection was done, as $\cardPprior$ shifts from $10,565$ $\left(\Nd = 8\right)$ to $1.7 \times 10^6$ $\left(\Nd = 40\right)$, demonstrating the 
efficiency of the subset selection step which activates only a small fraction of the dictionary. When $\Nd$ further increases from 40 to 100 for a given $\Nq$, the performance remains essentially the same with hardly any loss of accuracy: the solution method is able to capture the low-dimensional manifold onto which the solution essentially lies and an increase in the size of the solution space hardly affects the number of samples it requires. This capability is a crucial feature when available data are scarce and the solution space is very large.
As an illustration, when $\Nd = 100$, and with the parameters retained, the potential cardinality of the approximation basis is about $27 \times 10^6$ while the number of available samples is $\mathcal{O}\left(100-10,000\right)$. It clearly illustrates the pivotal importance of the subset selection step. Note that if one substitutes a PC approximation to the present HDMR format, about $352 \times 10^9$ terms need be evaluated with the present settings, a clearly daunting task.

For sake of completeness, the approximation given by a CP-format, Eq.~\eqref{CP-like}, is also considered. The univariate functions $\left\{f_{i,r}\right\}$ are approximated with the same polynomial approximation as in the present gLARS-ALS approach and a Tikhonov-based regularized ALS technique is used to determine each $f_{i,r}$ in turn given the others. Upon convergence, the next set of modes $\left\{f_{1,r+1}, \ldots, f_{\Nd,r+1}\right\}$ is evaluated until a maximum rank $\nr$ set by cross-validation. At each rank $r$, the best approximation, as estimated by cross-validation, is retained from a set of initial conditions and regularization parameter values. As can be appreciated from Fig.~\ref{L2_MM_schemes_Nds}, the number of samples required for a given approximation error is significantly larger than with the present gLARS-HDMR method.

\subsubsection{Subset selection}

\begin{figure}[!ht]
\begin{center}
  \includegraphics[angle = 0, width = 0.49\textwidth, draft = false]{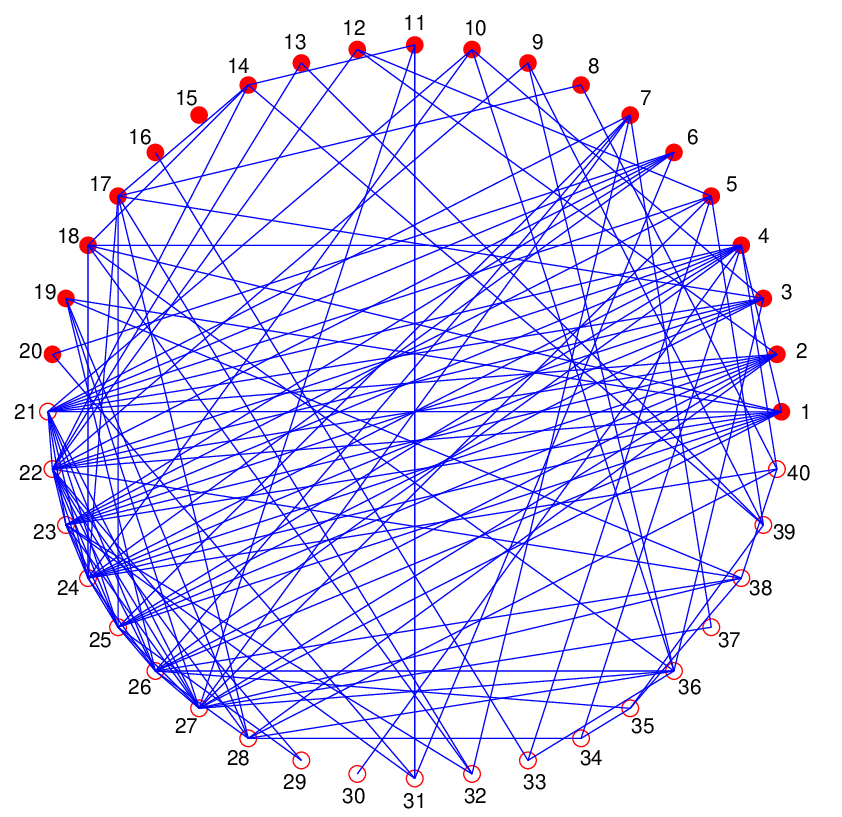}
  \includegraphics[angle = 0, width = 0.49\textwidth, draft = false]{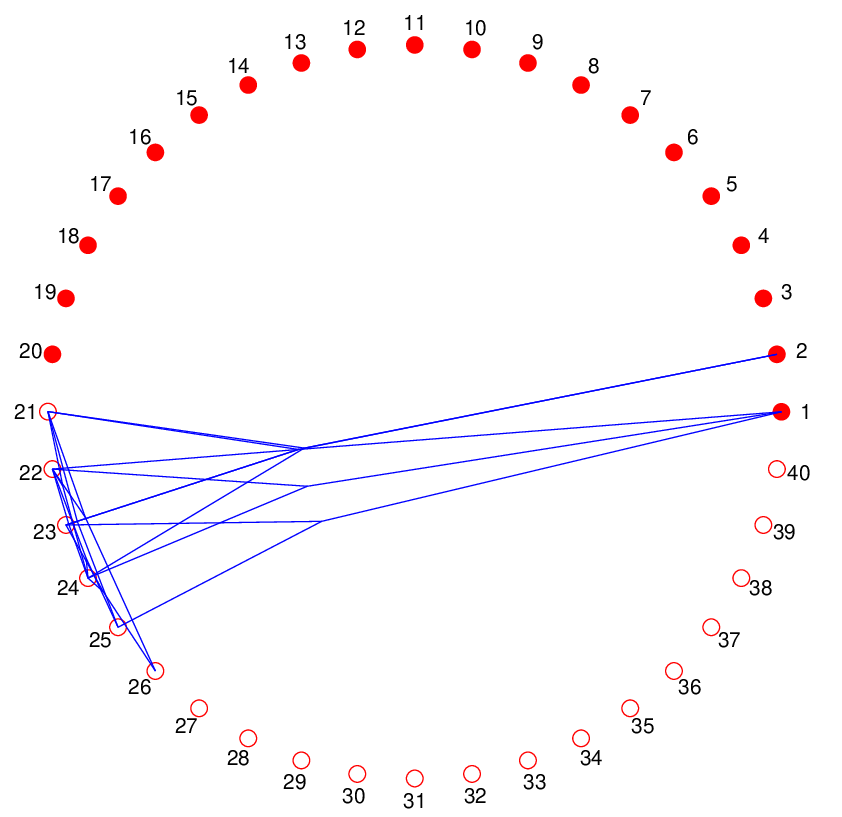}
  \caption{Graphical representation of the interaction modes retained by the subset selection procedure. Left: second order modes are plotted as a line linking two dimensions (bullets). Right: third order modes are represented as 3-branch stars and connect three dimensions.}
  \label{duo_tree}
\end{center}
\end{figure}

To further illustrate the subset selection step, the set of second and third order interaction retained modes $\left\{f_{\bgamma}\right\}_{\bgamma \in \basissetfpost}$ are plotted in Fig.~\ref{duo_tree} in the $\Nd = 40$ case. Each bullet represents one of the $\Nd$ stochastic dimensions and each line connects two (2-nd order, left plot) or three (3-rd order, right plot) dimensions, denoting a retained mode. The first $\Nd_F = 20$ of the 40 dimensions are associated with the source term $F$ in the stochastic equation and are represented as the solid bullets of the first two quadrants, $d \in \left[1, 20\right]$. The other $\Nd_\nu = 20$ dimensions are associated with the uncertain diffusion coefficient $\nu$ and are plotted as open bullets in the 3-rd and 4-th quadrants, $d \in \left[21, 40\right]$. The dimensions introduced by these two quantities are sorted with the associated magnitude of the eigenvalues $\sigma_F$ and $\sigma_\nu$ of their kernel, see Eqs.~\eqref{KL_eigen}, which decreases along the 
counter-clockwise direction. Hence, the norm of the eigenvalues of the kernel associated with $F$ decreases when one goes counter-clockwise from the first to the second quadrant. Likewise, the norm of the eigenvalues associated with dimensions introduced by $\nu$ decreases from the third to the fourth quadrant. Dominant dimensions of the stochastic space for the output $u$ approximation are thus expected to lie at the beginning of the first and/or third quadrant.

From the plot of second order modes (left), the subset selection process is seen to retain interaction modes mainly associated with dominant eigenvalues of both $F$ and $\nu$: they mainly link bullets from the first (dominant) dimensions associated with $F$ to the first (dominant) dimensions associated with $\nu$, as one might expect. Further, modes associated with two dimensions both introduced by $\nu$ are seen to be selected while two dimensions both associated with $F$ are rarely connected: the subset selection procedure is able to capture the nonlinearity associated with $\nu$ in the QoI and retains corresponding interaction modes. Indeed, note from Eq.~\eqref{diff_eq} that the source term $F$ interacts linearly with the solution $u$ while the diffusion coefficient is nonlinearly coupled with $u$ and hence, interaction modes between two dimensions introduced by $F$ do not contribute to the approximation. The third order modes (right plot) also illustrate the nonlinearity associated with $\nu$: the 
retained modes either connect dimensions associated with $\nu$ only or with one $F$-related and two $\nu$-related dimensions. Again, no two dimensions of $F$ are connected, consistently with the linear dependence of $u$ with $F$. These results illustrate the effectiveness of the procedure to unveil the dominant dependence structure and to discard unnecessary approximation basis functions.

%
%

\subsubsection{Robustness} \label{Robust_results}

The robustness of the approximation against measurement noise is now investigated. The dataset is corrupted with noise.
Denoting the nominal value with a star as superscript, noise in the coordinates is modeled as
\be
\bxiq = {\bxiq}^\star + s \, \bzetaq, \qquad \forall \, 1 \le q \le \Nq, \qquad s > 0.
\ee

The noise is modeled as an additive $\Nd$-dimensional, zero-centered, unit variance, Gaussian random vector $\bzeta$ biased so that $\bxiq \in \left[-1, 1\right]^\Nd$, $\forall \, q$. It is independent from one sample $q$ to another. Without loss of generality, measurements are here modeled as being corrupted with a multiplicative noise: $\uq = {\uq}^\star \, \left(1  + s_u \, \zeta_u^{(q)}\right)$, with $s_u = 0.2$ and $\zeta_u \sim \mathcal{N}\left(0, 1\right)$.

\begin{figure}[!ht]
\begin{center}
  \includegraphics[angle = -90, width = 0.5\textwidth, draft = false]{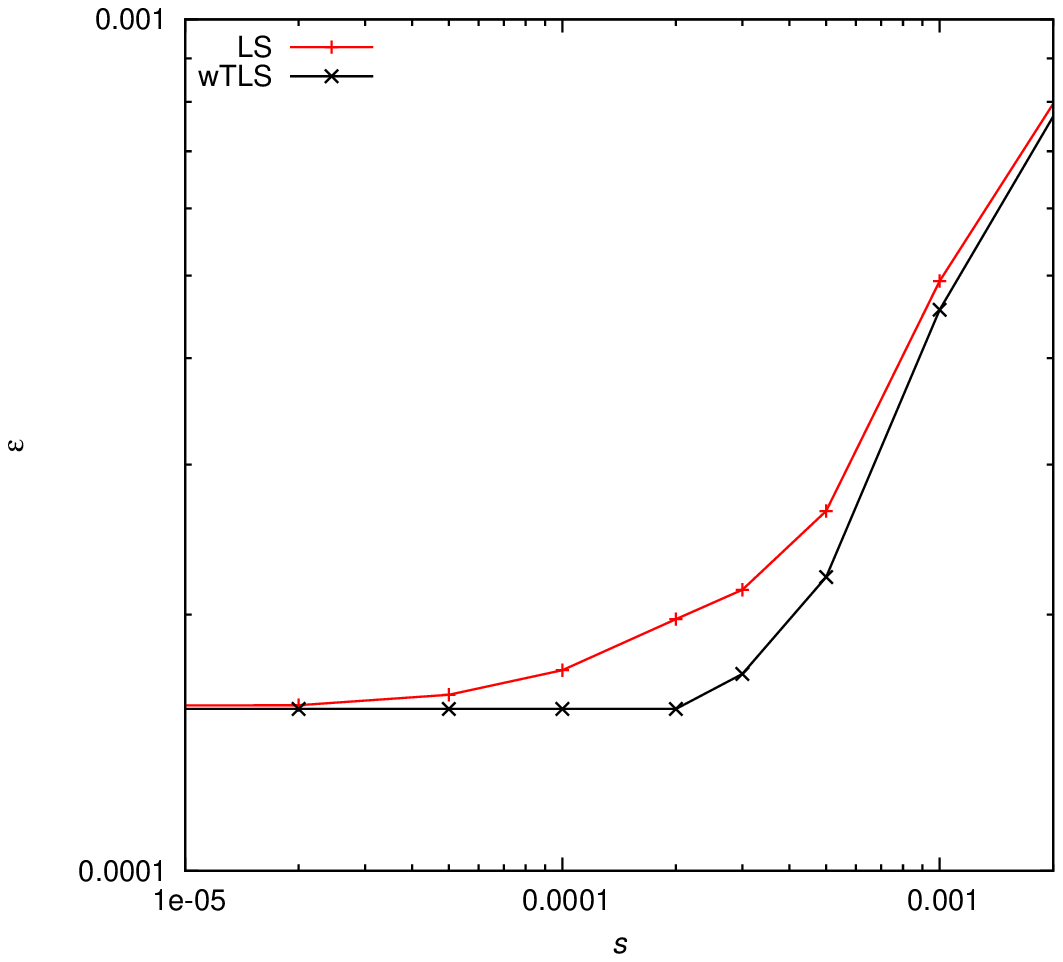}
  \caption{Robustness of the approximation w.r.t. noise in the data: approximation error $\varepsilon$ from the standard least squares (LS) and weighted total least squares (wTLS). $\Nd = 5$, $\Nq = 500$.}
  \label{Noise_fig}
\end{center}
\end{figure}

The evolution of the approximation accuracy when the noise intensity $s$ in the coordinates varies is plotted in Fig.~\ref{Noise_fig} in terms of error estimation $\varepsilon$. We compare gLARS-ALS using standard least squares (LS) with its `robust' counterpart relying on weighted total least squares (wTLS) as discussed in section~\ref{Noise_section}.

When the noise intensity increases, the error exponentially increases, quickly deteriorating the quality of the approximation with a noise standard deviation here as low as $s = 3 \times 10^{-5}$. When the noise is strong (low SNR), both the LS and the wTLS methods achieve poor accuracy. However, if the dataset is only mildly corrupted with noise, the wTLS approach is seen to achieve a significantly better accuracy than the standard least squares, while the solution process is significantly slower than that using the standard least squares. The present paper is based on the assumption that the critical part of the whole solution chain of determining a good approximation of the QoI is the data acquisition and that the cost of the post-processing part is not the main issue. However, while it is useful only on a range of SNR and somehow computationally costly, this feature is deemed important for a successful solution method in an experimental context where noise is naturally present.

If the noisy dataset is unbiased, possible further improvements upon the wTLS approach include lowering the complexity of the approximation model. Indeed, the well known bias-variance tradeoff indicates that a more robust, while less accurate, approximation can be obtained when the complexity of the retained model decreases. To improve the robustness of our present approach, a natural way is hence to trade some accuracy for some additional robustness. For instance, a predictor selection within each retained groups $\left\{f_{\bgamma \in \basissetfpost}\right\}$ can be considered, further lowering the final number of coefficients involved in the approximation and likely improving its robustness w.r.t. noise in the data. This could be achieved by estimating the approximation coefficients via a \emph{penalized} (total) least squares problem as mentioned in section~\ref{subsection_regression}.

\subsubsection{Scaling of the solution algorithm} \label{Num_cost_example}

In this section, the numerical complexity associated with the different steps of the solution method is illustrated in terms of computational time. Numerical experiments are carried-out with varying number of samples $\Nq$ and solution space dimensions $\Nd$. When one is varying, the other remains constant. The nominal parameters are $\Nd = 40$ (dimension of the stochastic space), $\No = 6$ (maximum total order of the Legendre polynomials), $\Nq = 1000$ (number of samples), $\Ninter = 3$ (maximum interaction order of the truncated HDMR approximation), $\NoLARS = 5$ (maximum total polynomial order in the subset selection step).

\begin{figure}[!ht]
\begin{center}
  \includegraphics[angle = -90, width = 0.49\textwidth, draft = false]{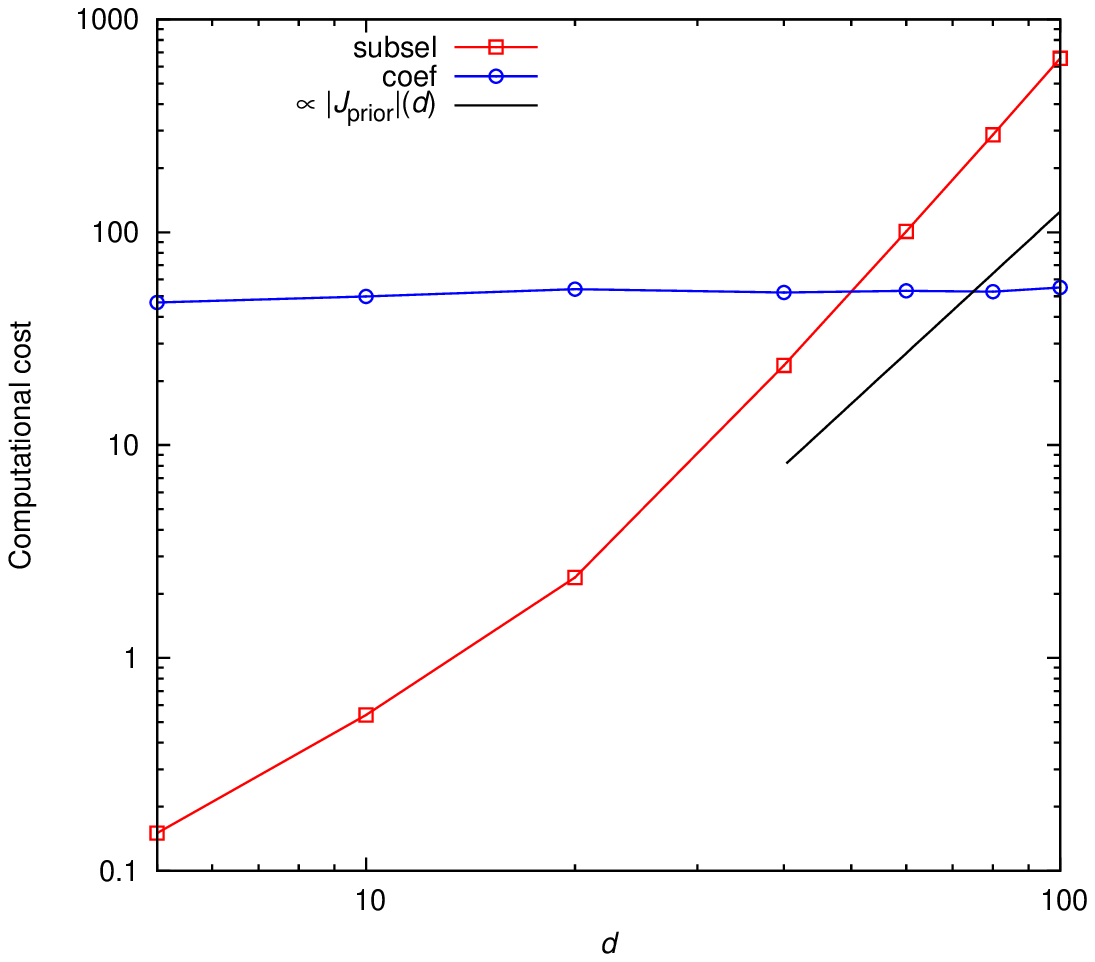}
  \includegraphics[angle = -90, width = 0.49\textwidth, draft = false]{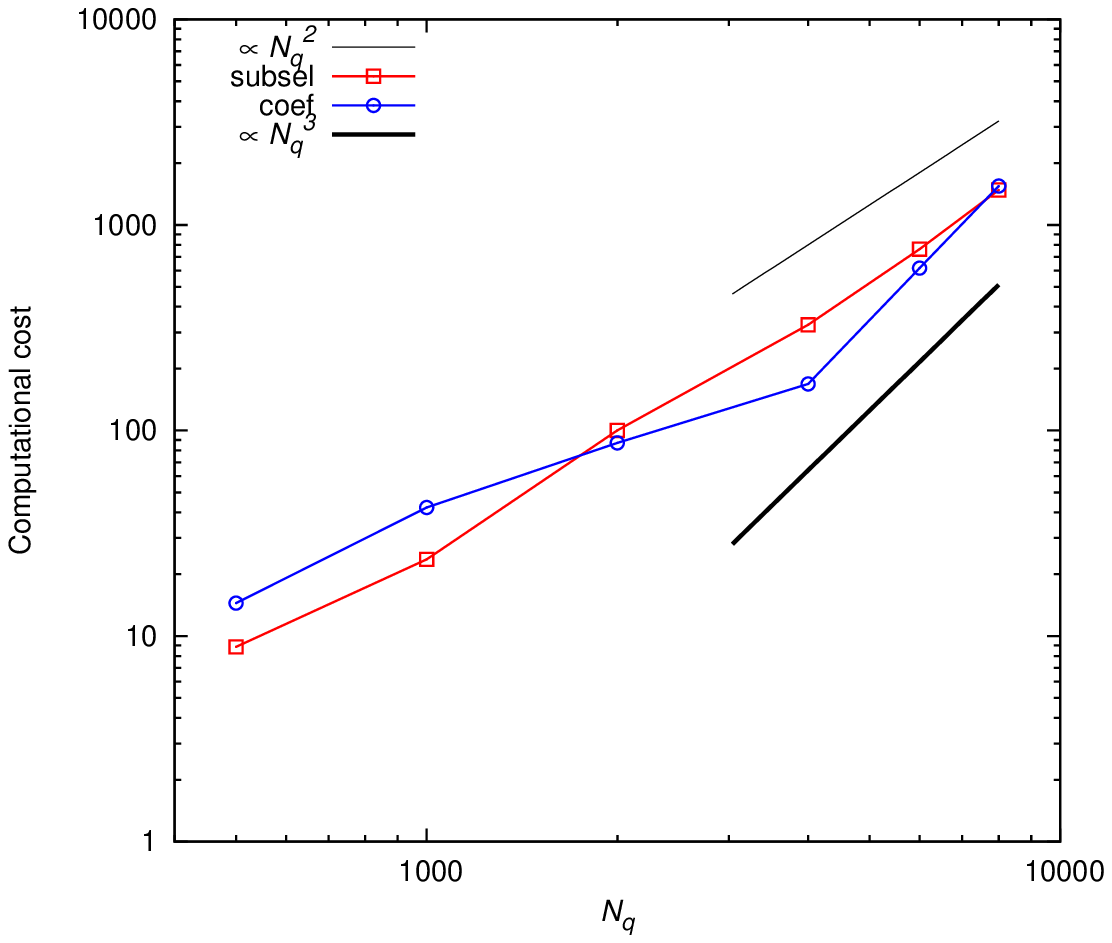}
  \caption{Numerical cost of the subset selection and coefficients evaluation steps as a function the stochastic dimension $\Nd$ and size of the dataset $\Nq$. Approximation coefficients are fully updated for each new mode. Nominal parameters are $\Nd = 40$, $\No = 6$, $\Nq = 1000$, $\Ninter = 3$, $\NoLARS = 5$, $\NinterPC = 3$.}
  \label{Num_cost}
\end{center}
\end{figure}

Numerical results are gathered in Fig.~\ref{Num_cost}.
The asymptotic behavior of the number $\nJ$ of required subset selection iterations as introduced in section \ref{Complexity_sec} might be different according to which limit is considered. For the present stochastic diffusion problem, first and second interaction order modes tend to be selected first. Assuming the active set $\basissetfpost$ is dominated by first and second interaction order modes, it can easily be shown that the number of retained groups then satisfies
\be
\nJ \le 1 + \nJzero + \textrm{min}\left[\frac{\Nd \, \left(\Nd-1\right)}{2}, 2 \, \frac{\Nq - \nJzero \, \NoLARS}{\NoLARS \, \left(\NoLARS+1\right)}\right], \qquad \nJzero \le \textrm{min}\left[\Nd, \frac{\Nq}{\NoLARS}\right].
\ee

In the present example, second order interaction groups dominate the retained set so that the number of retained groups tends to scale as $\nJ \propto \Nq / \NoLARS^2$.
From Eq.~\eqref{theo_cost_subset} and for the present nominal parameters, it results in the following limit behavior for the subset selection step:
\bea
\lim_{\Nd \rightarrow +\infty} \mathscr{J}_\mathrm{subsel} & \propto & \Nq^2 \, \cardPprior / \NoLARS^2 \qquad \longrightarrow {\rm here:} \quad \propto \Nq^2 \, \Nd^\Ninter \, \NoLARS^{\Ninter - 2}, \nonumber \\
\lim_{\Nq \rightarrow +\infty} \mathscr{J}_\mathrm{subsel} & \propto & \Nq^2 \, \cardPprior / \NoLARS^2 \qquad \longrightarrow {\rm here:} \quad \propto \Nq^2 \, \Nd^\Ninter \, \NoLARS^{\Ninter - 2}.
\eea

Similarly, the cost associated with the coefficients evaluation is considered. The number of interaction modes $\basissetfcur$ effectively varies between $1$ and $\mathcal{O}\left(\Nq / \No^2\right)$ along the solution procedure, and, since the cost associated with solving the least squares problem dominates that of the matrix assembly, the cost of their evaluation finally simplifies in $\mathscr{J}_\mathrm{coef} \propto \mathcal{O}\left(\Nq^2 \, \No\right)$ or $\mathscr{J}_\mathrm{coef} \propto \mathcal{O}\left(\Nq^3 / \No\right)$ depending on whether the coefficients are updated whenever an additional group is considered or not, see section~\ref{Algo_lambda} and step~(\ref{Updateornot}) in Algorithm~\ref{Algo_sumup}. In the present regime, the cost is found not to depend on $\Nd$.

These asymptotic behaviors are consistent with the numerical experiments as can be appreciated from Fig.~\ref{Num_cost}. The coefficients are here updated whenever a new mode from the selected set is considered, hence $\mathscr{J}_\mathrm{coef} \propto \mathcal{O}\left(\Nq^3 / \No\right)$.
It is seen that the subset selection step scales less favorably than the coefficients evaluation step with the dimensionality of the random variable. This stresses the benefit of a carefully chosen {\apriori} approximation basis to reduce as much as possible the cardinality $\cardPprior$.

%
%

\subsection{Approximation of the solution random field}

We now consider the approximation of the space-dependent random solution $u\left(x, \bxi\right)$ under the form \eqref{KL_approx_format} using Algorithm~\ref{Algo_separated}. The approximation obtained from different number of samples $\left\{\xq, \bxiq, \uq\right\}$ is compared with the Karhunen-Lo\`eve modes, computed from a full knowledge of the QoI, hereafter referred to as the reference solution.\footnote{The spatial $\left\{\modex_n\right\}$ and stochastic modes $\left\{\lambda_n\right\}$ are sequentially determined from \eqref{x_modes} via an ALS approach. Since the decomposition is two-dimensional, $u\left(\bx, \bxi\right) \approx \sum_{n=0}^{\KLrank} { \modex_n\left(\bx\right) \, \modesto_n\left(\bxi\right) }$, the approximation problem is convex, see for instance \cite{Grasedyck_10}, and the ALS approach converges to the best rank-1 approximation of the matricized $\bu$ in the Frobenius sense. If the data-driven inner product $\left< \cdot, \cdot \right>_{\Nq}$ was inducing a cross-norm (it only induces a semi-norm), then $\left< \modex \, \modesto, \modex \, \modesto \right>_{\Nq} = \left\|\bmodex\right\|_2^2 \, \left\|\bmodesto\right\|_2^2$ and the pair $\left(\bmodex, \bmodesto\right)$ would be the dominant rank-1 approximation of the matricized $\bu$. The Karhunen-Lo\`eve decomposition of $u$ is thus the reference solution one should obtain in the particular case where the empirical inner product induces a cross-norm and $\Nq \rightarrow \infty$.}

\begin{figure}[!ht]
\begin{center}
  \includegraphics[angle = -90, width = 0.49\textwidth, draft = false]{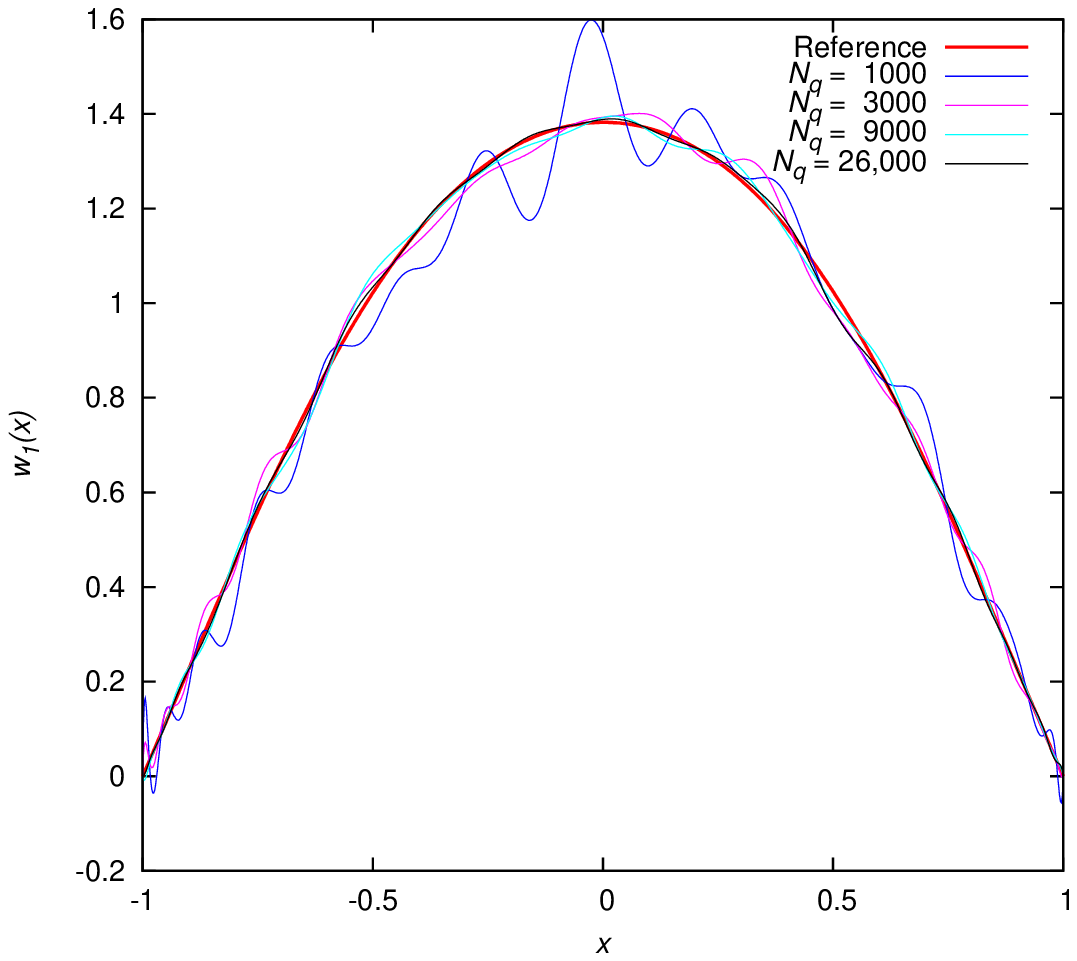}
  \includegraphics[angle = -90, width = 0.49\textwidth, draft = false]{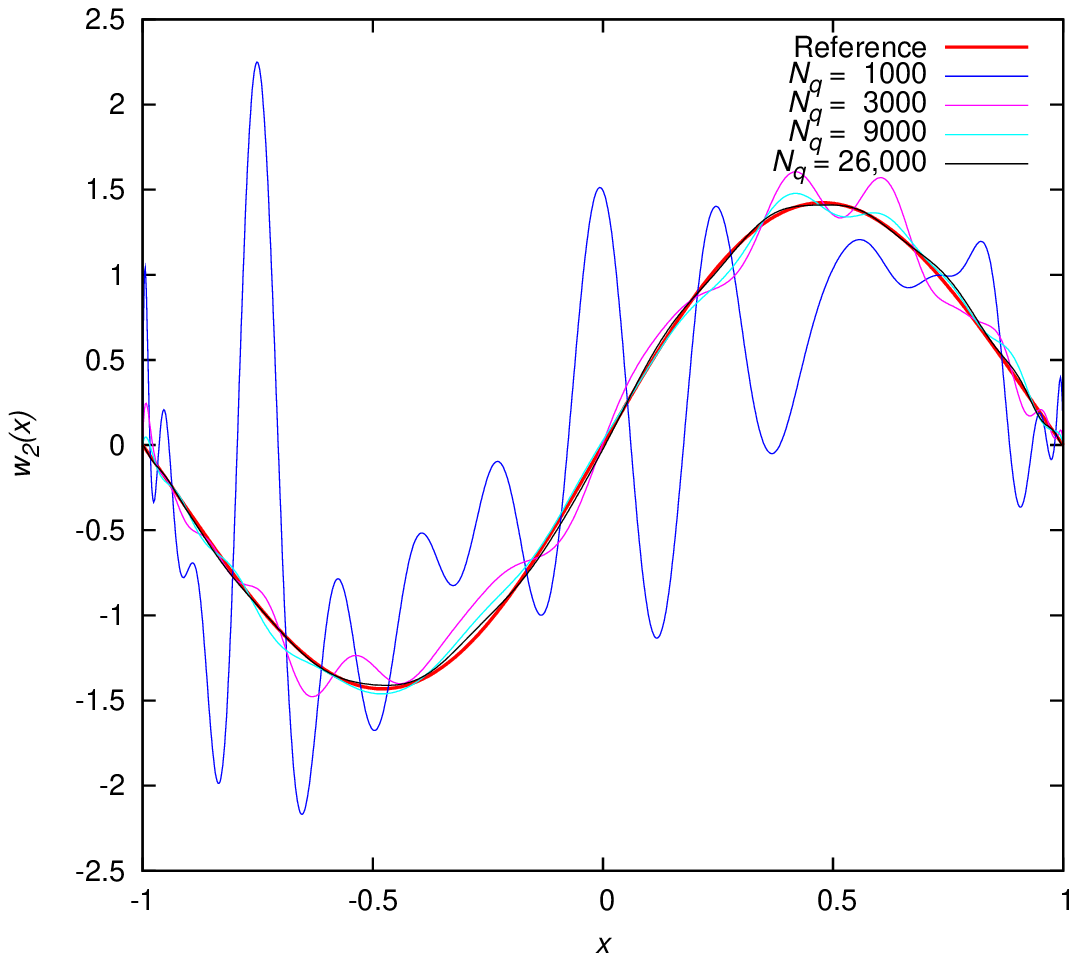}
  \caption{First ($\modex_1\left(x\right)$, left) and second ($\modex_2\left(x\right)$, right) spatial approximation modes of the stochastic diffusion solution. The reference (Karhunen-Lo\`eve) solution is plotted for comparison (thick line).}
  \label{KL_1D}
\end{center}
\end{figure}

The simulation relies on the following parameters: $\cardx = 32$, $\No = 10$, $\Nd = 6$, $\Ninter = 3$, $\NinterPC = 3$. The potential approximation basis cardinality is about $\cardx \, \cardPprior \simeq 10^5$. Fig.~\ref{KL_1D} shows the first and second spatial modes, $\modex_1\left(x\right)$ and $\modex_2\left(x\right)$ for different sizes of the dataset, $\Nq = 1000$, $3000$, $9000$ and $26,000$. The mean mode $\modex_0\left(x\right)$ is virtually indistinguishable from the reference solution mean mode for any of the dataset sizes and is not plotted. On the left plot $\left(\modex_1\left(x\right)\right)$, it is seen that the approximation is decent, even with as low as $\Nq = 1000$ samples. For $\Nq = 3000$, the approximation is good. This $\left(1 + \Nd\right) = 7$-dimensional case corresponds to $\Nq^{1 / (1 + 6)} \simeq 3.1$ samples per solution space dimension only and about $\Nq / \left(\cardx \, \cardPprior\right) \simeq 3 \%$ of the potentially required information.

For approximating the second spatial mode (Fig.~\ref{KL_1D}, right plot), more points are needed to reach a good accuracy but $\Nq = 26,000$ is seen to already deliver a good performance. Quantitative approximation error results are gathered in Table~\ref{epstab_1D} for various separation ranks $\KLrank$ and number of samples $\Nq$.

\begin{table}[ht!]
\begin{center}
\begin{tabular}{r|ccc}
\hline \hline
$\Nq$ $\backslash$ $\KLrank$ & 0 & 1       & 2       \\
\hline
1000  & $5.5 \times 10^{-3}$ & $7.4 \times 10^{-4}$ & $7.4 \times 10^{-4}$ \\
3000  & $5.5 \times 10^{-3}$ & $4.2 \times 10^{-4}$ & $2.7 \times 10^{-4}$ \\
9000  & $5.5 \times 10^{-3}$ & $3.1 \times 10^{-4}$ & $1.0 \times 10^{-4}$ \\
26,000 & $5.4 \times 10^{-3}$ & $2.8 \times 10^{-4}$ & $6.2 \times 10^{-5}$ \\
\hline
\end{tabular}
\caption{Evolution of the approximation error $\varepsilon$, as defined in Eq.~\eqref{sto_error_def}, with the decomposition rank $\KLrank$ and the number of samples $\Nq$.}
\label{epstab_1D}
\end{center}
\end{table}

The satisfactory performance of the present method can be understood from the upper part of the Karhunen-Lo\`eve approximation (normalized) spectrum plotted in Table~\ref{KLspectrum_1D}. The norm of the eigenvalues decays quickly so that the first two modes contribute more than 90 \% of the QoI $L^2$-norm, showing that this problem efficiently lends itself to the present separation of variables-based methodology.

\begin{table}[ht!]
\begin{center}
\begin{tabular}{l|cccccccccc}
\hline \hline
$i$ & 1 & 2 & 3 & 4 & 5 & 6 & 7 & 8 & 9 & 10 \\
\hline
$\sigma_i$ & 144 & 30.2 & 13.6 & 2.64 & 1.37 & 0.250 & 0.0817 & 0.0167 & 0.00891 & 0.00232 \\
\hline
\end{tabular}
\caption{Normalized upper spectrum of the Karhunen-Lo\`eve approximation.}
\label{KLspectrum_1D}
\end{center}
\end{table}

\subsection{A Shallow Water flow example}

The methodology is now applied to the approximation of the stochastic solution of a Shallow Water flow simulation with multiple sources of uncertainty. It is a simple model for the simulation of wave propagation on the ocean surface. Waves are here produced by the sudden displacement of the sea bottom at a given magnitude in time, extension and location, all uncertain.

\subsubsection{Model}

The problem is governed by the following set of equations:
\bea
\frac{D \, \vx}{D \: t} & = & f_C \: \vy - g \: \frac{\partial h}{\partial x_1} - b \: \vx + \Su,  \label{SWE_ori_a}\\
\frac{D \: \vy}{D \: t} & = & - f_C \: \vx - g \: \frac{\partial h}{\partial x_2} - b \: \vy + \Sv, \\
\frac{\partial h}{\partial t} & = & - \frac{\partial \left( \vx \: \left( H + h \right) \right)}{\partial x_1} - \frac{\partial \left( \vy \: \left( H + h \right) \right)}{\partial x_2} + \Sh, \label{SWE_ori_c}
\eea
where $\left(\vx\left(\bx,\bxi,t\right) \, \vy\left(\bx,\bxi,t\right)\right)$ is the velocity vector at the surface, $\bx = \left(x_1 \, x_2\right) \in \Omega \subset \mathbb{R}^2$, $h\left(\bx,\bxi,t\right)$ the elevation of the surface from its position at rest, $H\left(\bx\right)$ the sea depth, $f_C$ models the Coriolis force, $b$ is the viscous drag coefficient, $g$ the gravity constant and $\Su\left(\bx,\bxi,t\right)$, $\Sv\left(\bx,\bxi,t\right)$, $\Sh\left(\bx,\bxi,t\right)$ are the source fields. Without loss of generality, the drag $b$ and the Coriolis force $f_C$ are neglected. No slip boundary conditions apply for the velocity. The sources are modeled as acting on $h$ only, $\Su \equiv 0$ and $\Sv \equiv 0$. $\Sh$ models the source term acting on $h$ due to, say, an underwater seismic event. The fluid density and the free surface pressure are implicitly assumed constant. Full details on the numerical implementation of a similar problem are given in \cite{Mathelin_al_CompMech}.

\subsubsection{Sources of uncertainty}

Let $\bxi = \left( \bxi' \, \bxi''\right)$. The source $\Sh$ is uncertain and is modeled as a time-dependent, spatially distributed, quantity:
\be
\Sh\left(\bx,\bxi,t\right) = a_t\left(\bxi',t\right) \, a_\xi\left(\xi''_1\right) \, \exp\left({-\frac{\left(\bx - \bx_\Sh\left(\xi''_3\right)\right)^T \, \left(\bx - \bx_\Sh\left(\xi''_3\right)\right)}{\sigma_\Sh\left(\xi''_2\right)^2}}\right),
\ee
where $a_t\left(\bxi',t\right)$ is a given time envelop, $a_\xi\left(\xi''_1\right)$ the uncertain source magnitude, $\sigma_\Sh\left(\xi''_2\right)$ drives the uncertain source spatial extension and $\bx_\Sh\left(\xi''_3\right)$ is the uncertain spatial location. The time envelop $a_t\left(\bxi', t\right)$ is described with a $\NH$-term expansion:
\be
a_t\left(\bxi',t\right) = \overline{a_t}\left(t\right) + \sum_{i=1}^{\NH}{\sqrt{\lambda_i} \: \xi_i' \left(\theta\right) \, \varphi_i^{a_t}(t)}, \label{H_eq}
\ee
with $\bxi' = \left(\xi'_1 \, \ldots \, \xi'_\NH\right)$ the stochastic germ associated to the uncertainty in $a_t$. Random variables $\left\{\xi'_i\right\}_{i=1}^{\NH}$ are \textit{iid}, uniformly distributed. The solution of the Shallow Water problem then lies in a $\left(\Nd = \NH + 3\right)$-dimensional stochastic space.

\subsubsection{Approximation from an available database}

As an illustration of the methodology, we aim at approximating the sea surface field at a fixed amount of time $t^\star$ after a seismic event. The QoI is then a random field $u\left(\bx, \bxi\right) = h\left(\bx, \bxi, t^\star\right)$. An accurate description of this field is of importance for emergency plans in case of a seaquake. Sea level measurements of the surface at various spatial locations from past events constitute the dataset $\left\{\bxq, \bxiq, h\left(\bxq, \bxiq, t^\star\right)\right\}_{q=1}^{\Nq}$ used to derive an approximation of $u$ under a separated form: $u\left(\bx, \bxi\right) \approx \left<u\right>_{\Nq}\left(\bx\right) + \sum_{n=1}^\KLrank{\modex_n\left(\bx\right) \, \modesto_n\left(\bxi\right)}$.

The solution method here relies on a $\Nq = 37,000$-sample dataset complemented with $\Nqhat = 5000$ cross-validation samples and a $\Nqt = 5000$ set for error estimation. We consider a $\NH=5$ expansion for the time envelop, leading to a stochastic dimension of $\Nd = 5+3=8$. The effective number of samples per dimension is then about $\Nq^{1/\left(\Ndx + \Nd\right)} \simeq 2.9$. The approximation is determined based on a $\cardx = 484$ spatial discretization DOFs (spectral elements) at the deterministic level and $\No=6$-th order Legendre polynomials $\left\{\psi_\alpha\right\}$, $\Ninter = 3$, $\NinterPC = 3$, for the stochastic modes. The cardinality of this {\apriori} basis is then $\cardx \, \cardPprior \simeq 770 \times 10^3 \gg \Nq$, again relying on an efficient subset selection step to make the approximation problem well-posed.

\begin{figure}[!ht]
\begin{center}
  \includegraphics[angle = -90, width = 0.49\textwidth, draft = false]{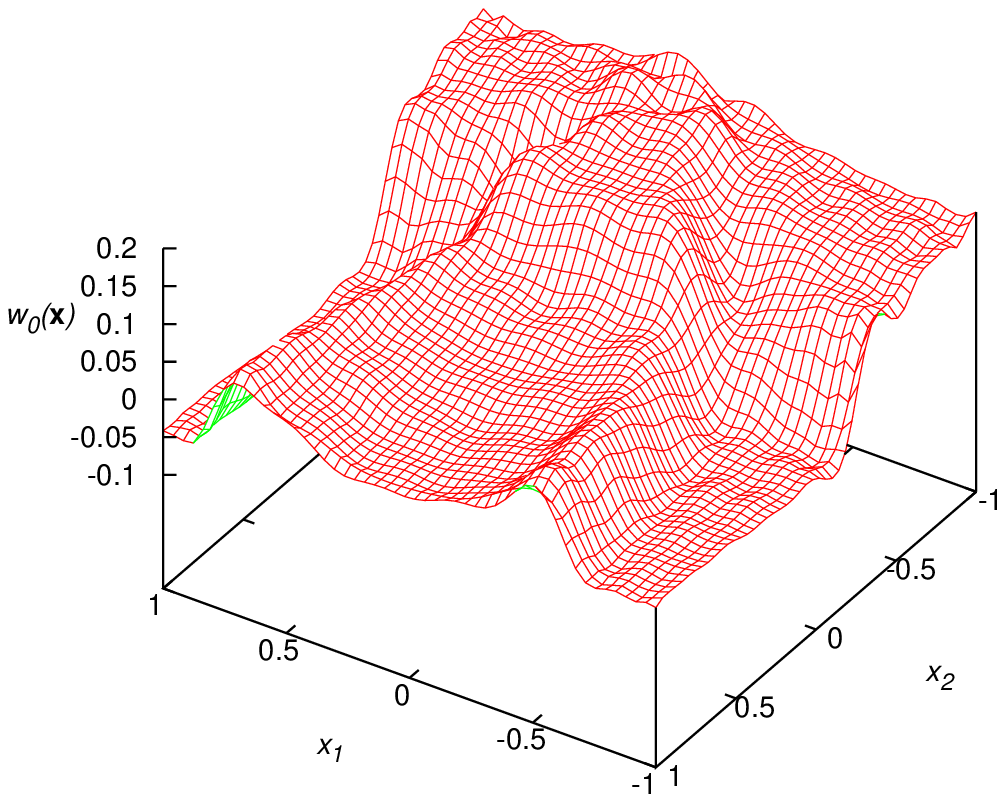}
  \includegraphics[angle = -90, width = 0.49\textwidth, draft = false]{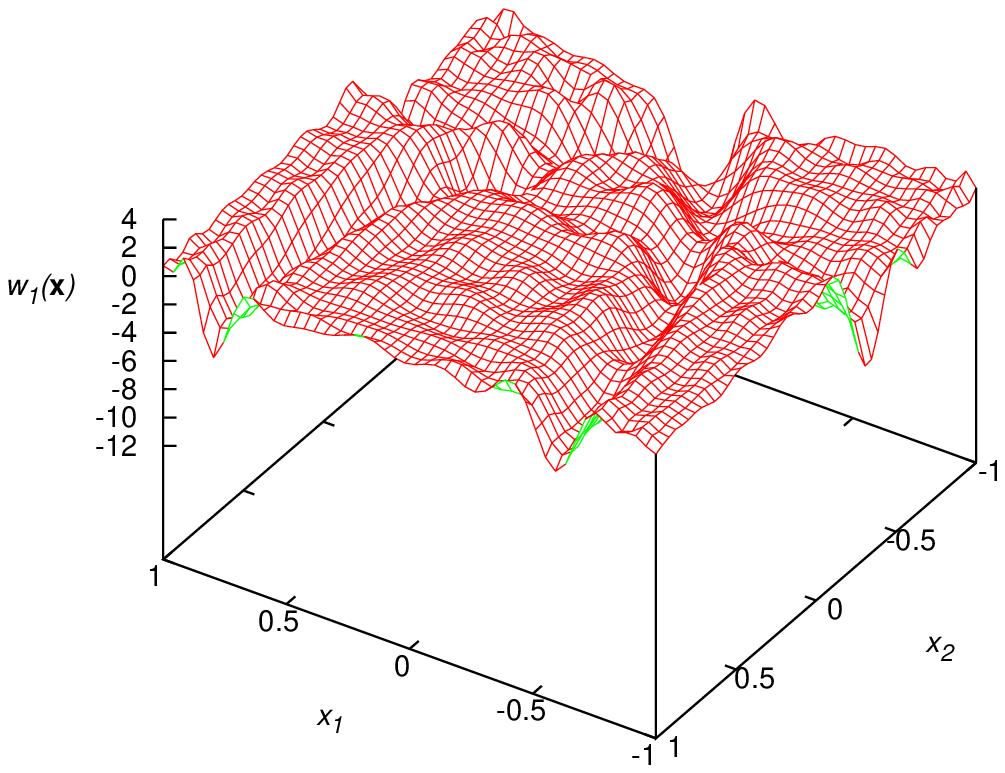}
  \caption{Mean ($\modex_0\left(\bx\right) \equiv \left<u\right>_{\Nq}\left(\bx\right)$, left) and first ($\modex_1\left(\bx\right)$, right) spatial modes.}
  \label{KL_SWE}
\end{center}
\end{figure}

The approximation error when the rank $\KLrank$ varies is shown in Table~\ref{Norme_x_SWE}. It is seen that estimating the mean spatial mode $\modex_0$ leads to a relative error of about $0.12$ while adding the first $\left(\modex_1, \modesto_1\right)$ and second $\left(\modex_2, \modesto_2\right)$ pair drops it to about $0.05$. Further adding pairs does not lower the approximation error with this dataset and more samples are needed to accurately estimate them. Spatial modes $\modex_0$ and $\modex_1$ of the separated approximation are plotted in Fig.~\ref{KL_SWE} for illustration.

\begin{table}[ht!]
\begin{center}
\begin{tabular}{c|cccc}
\hline \hline
$\KLrank$ & 0 & 1 & 2 & 3\\
\hline
$\varepsilon$ & 0.117 & 0.056 & 0.046 & 0.044 \\
\hline
\end{tabular}
\caption{Relative approximation error $\varepsilon$ evolution with the decomposition rank $\KLrank$. $\Nq = 37,000$.}
\label{Norme_x_SWE}
\end{center}
\end{table}

\section{Conclusion} \label{Conclusion_section}

In this paper, a methodology was proposed for deriving a functional representation of a random process only known through a collection of its pointwise evaluations. The proposed method essentially relies on an efficient determination of an approximation basis consistent with the available information. This involves the choice of an {\apriori} canonical HDMR format combined with tuning the basis \textit{via} a data-driven subset selection step. This subset selection is carried-out in a bottom-to-top manner, as opposed to a top-to-bottom manner as done in the Compressed Sensing standard framework. It essentially sorts the HDMR modes (groups of predictors) by their contribution in approximating the Quantity of Interest. The final approximation can rely on a different functional description of the modes, typically of higher order and/or nonlinear in the coefficients.

The method is progressive, data-driven, and was shown to here outperform current approximation techniques in terms of accuracy for a given number of samples. Its efficiency was demonstrated on two examples which have shown its ability to achieve a good approximation accuracy from a small dataset, as long as the quantity at hand is essentially lying on a low-dimensional manifold. In particular, the dominant dimensions are naturally revealed so that all the available information can be dedicated to approximate relevant dependences only.
Through a total least squares approach, it was also shown that some robustness can be achieved, an important feature if the dataset comes from experiments. Using a robust approximation was shown to bring up to a 2-fold improvement upon the approximation error using standard least squares, but at the price of a computational overhead.
The global solution method scales reasonably well, exhibiting a linear dependence with the cardinality of the {\apriori} basis dictionary and a quadratic or cubic dependence with the number of samples, depending on the coefficients update strategy.

The present work was focused on a general methodology, disregarding fine-tuning aspects. Among other things, a natural improvement would be to carry-out a predictor selection within each retained groups $\left\{f_{\bgamma \in \basissetfpost}\right\}$, further lowering the number of coefficients involved in the approximation. Moreover, the tensor structure of the Hilbert stochastic space can be exploited and developments towards a data-driven multilinear algebra effective tool for high-dimensional uncertainty quantification are currently carried-out.

\section*{Acknowledgement}
The author gratefully acknowledges Tarek El Moselhy and Faidra Stavropoulou for stimulating discussions and useful comments. This work is part of the TYCHE project (ANR-2010-BLAN-0904) supported by the French Research National Agency (ANR). 

\bibliographystyle{lionel_style}
\bibliography{biblio}

\appendix
\section{A motivating example} \label{Appendix_A}

To assess the choice of our {\apriori} functional form for approximating a random variable, and while choosing a good basis is problem-dependent, let us consider a simple motivating example in the form of the 1-D stochastic diffusion equation presented in section~\ref{StoDifEq}, briefly recalled here for sake of convenience:
\be
\nabla_x \, \left(\nu\left(x, \bxi\right) \, \nabla_x u\left(x,\bxi\right)\right) = F\left(x, \bxi\right), \qquad u\left(x_-,\bxi\right) = u_-, \: u\left(x_+,\bxi\right) = u_+. \label{diff_eq2}
\ee

The solution $u$ is approximated under a separated format $\displaystyle u\left(x, \bxi\right) \approx \sum_{n=0}^\KLrank{\modex_n\left(x\right) \, \modesto_n\left(\bxi\right)}$. The approximation space for the spatial modes $\left\{\modex_n\left(x\right)\right\}$ is given and we here focus on the accuracy of the approximation with different representations for the stochastic modes $\left\{\modesto_n\left(\bxi\right)\right\}$. Each stochastic mode is determined either in a CP-like format, Eq.~\eqref{CP-like}, or as a HDMR decomposition Eqs.~\eqref{HDMR_full}. In the latter case, interaction modes $\left\{f_\bgamma\right\}$ are approximated with a low-rank canonical decomposition on tensorized, unit-normed, univariate polynomials of maximum degree $\No$: $f_{\bgamma}\left(\left\{\xi_i\right\}_{i \in \bgamma}\right) \approx \sum_{r=1}^\nr \, \prod_{i \in \bgamma}{\sum_{\alpha=2}^\No{\coef_{\bgamma, \alpha}^{r, i} \, \psi_\alpha\left(\xi_i\right)}}$. This approximation is hereafter referred to as a CP-HDMR decomposition. Similarly, univariate functions $\left\{f_{i, r}\right\}_{i=1}^\Nd$ involved in the CP decomposition Eq.~\eqref{CP-like} are approximated with the same polynomials: $f_{i, r}\left(\xi_i\right) \approx \sum_{\alpha = 1}^{\No}{\coef_{\alpha, i, r} \, \psi_\alpha\left(\xi_i\right)}$.


The representation of the stochastic modes here relies on $\No = 8$-th order univariate Legendre polynomials $\left\{\psi_\alpha\right\}$. The dimension of the problem is chosen to be $\Nd_\nu = \Nd_F = 5$ so that $\Nd = 10$.

The CP-HDMR expansion is here built sequentially, starting with 0-th and 1-st order interaction modes only. From this first approximation of the output, the set of dominant dimensions is estimated from the $L^2$-norm of univariate interaction modes $\left\{f_\bgamma\right\}_{|\bgamma|=1}$. Only second order interaction modes $\left\{f_\bgamma\right\}_{|\bgamma|=2}$ in these dominant dimensions are next estimated and the set of dominant dimensions is then further refined based on both 1-st and 2-nd order interaction modes via the sensitivity Sobol indices, see \ref{Sobol_sec}. Third order modes are then computed for this new set of dominant dimensions only and the procedure is repeated until some stopping criterion is met, for instance a maximum interaction order $\Ninter$ or a maximum basis cardinality $\cardP$. The number of samples $\Nq$ is here chosen sufficiently large so that full knowledge on $u$ can be assumed. The approximation error then only comes from the choice of the approximation basis format, allowing a comparison.
This section is loose on details, focusing on the main conclusions and leaving more in-depth discussion for main text sections.

\begin{figure}[!ht]
\begin{center}
  \includegraphics[angle = -90, width = 0.5\textwidth, draft = false]{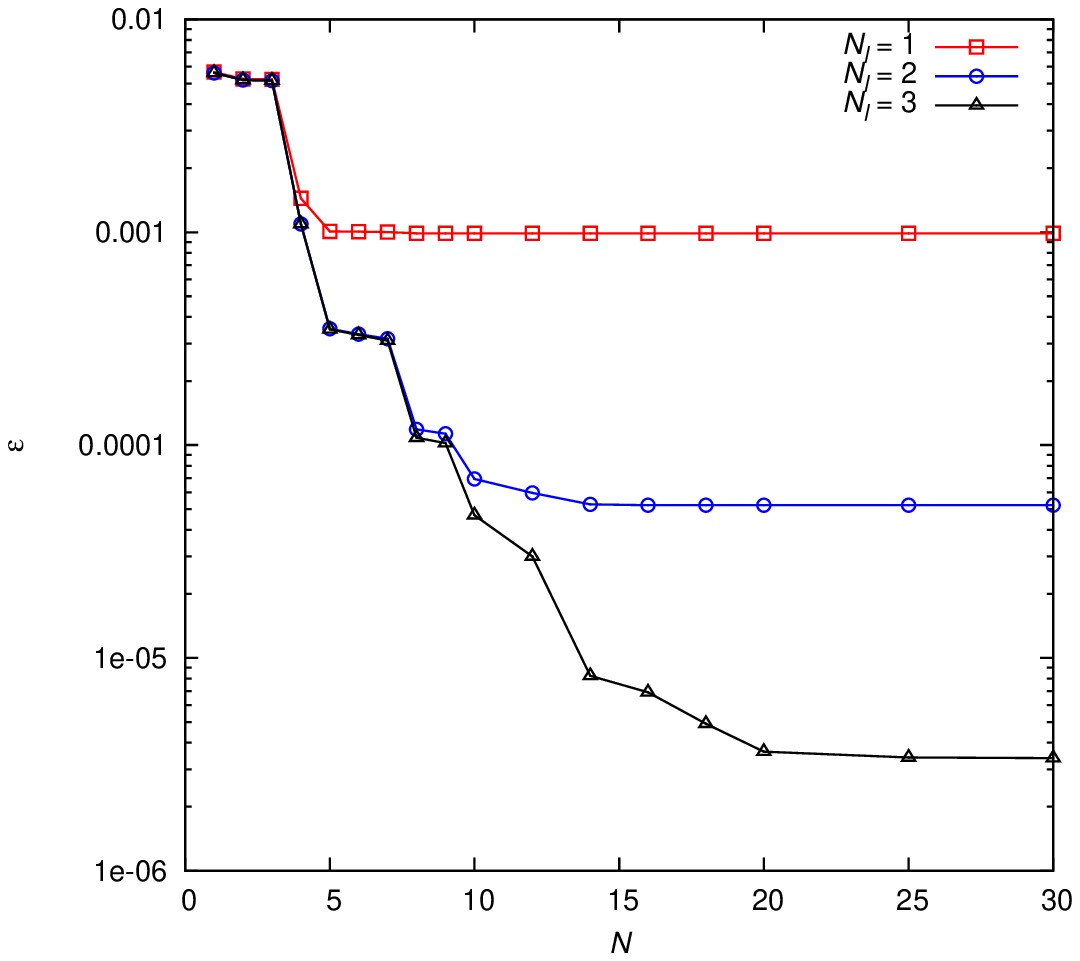}
  \caption{Convergence of the error estimation of the approximation with the order $\KLrank$ of the separated representation and the maximum interaction order $\Ninter$ of the HDMR expansion.}
  \label{L2_order_HDMR}
\end{center}
\end{figure}

First, the accuracy of the CP-HDMR approximation as a function of the decomposition rank $\KLrank$ is studied in terms of $\varepsilon$, Eq.~\eqref{sto_error_def}.
Plotted in Fig.~\ref{L2_order_HDMR}, the approximation error estimation $\varepsilon$ decreases when the maximum interaction order $\Ninter$ increases from 1 to 3 and as the decomposition rank $\KLrank$ increases.

The approximation is seen to improve exponentially fast as the number of modes $\KLrank$ in the separated representation increases until it reaches a plateau. Increasing the interaction order leads to an improved approximation: increasing from first to second order brings more than a one-order of magnitude improvement in the approximation error estimation and an additional order of magnitude from $\Ninter = 2$ to $\Ninter =3$. In this $\Nd = 10$ example, the approximation hence exhibits a high convergence rate with $\Ninter$, supporting our assumption that low-order interactions dominate the HDMR decomposition.

\begin{figure}[!ht]
\begin{center}
  \includegraphics[angle = -90, width = 0.49\textwidth, draft = false]{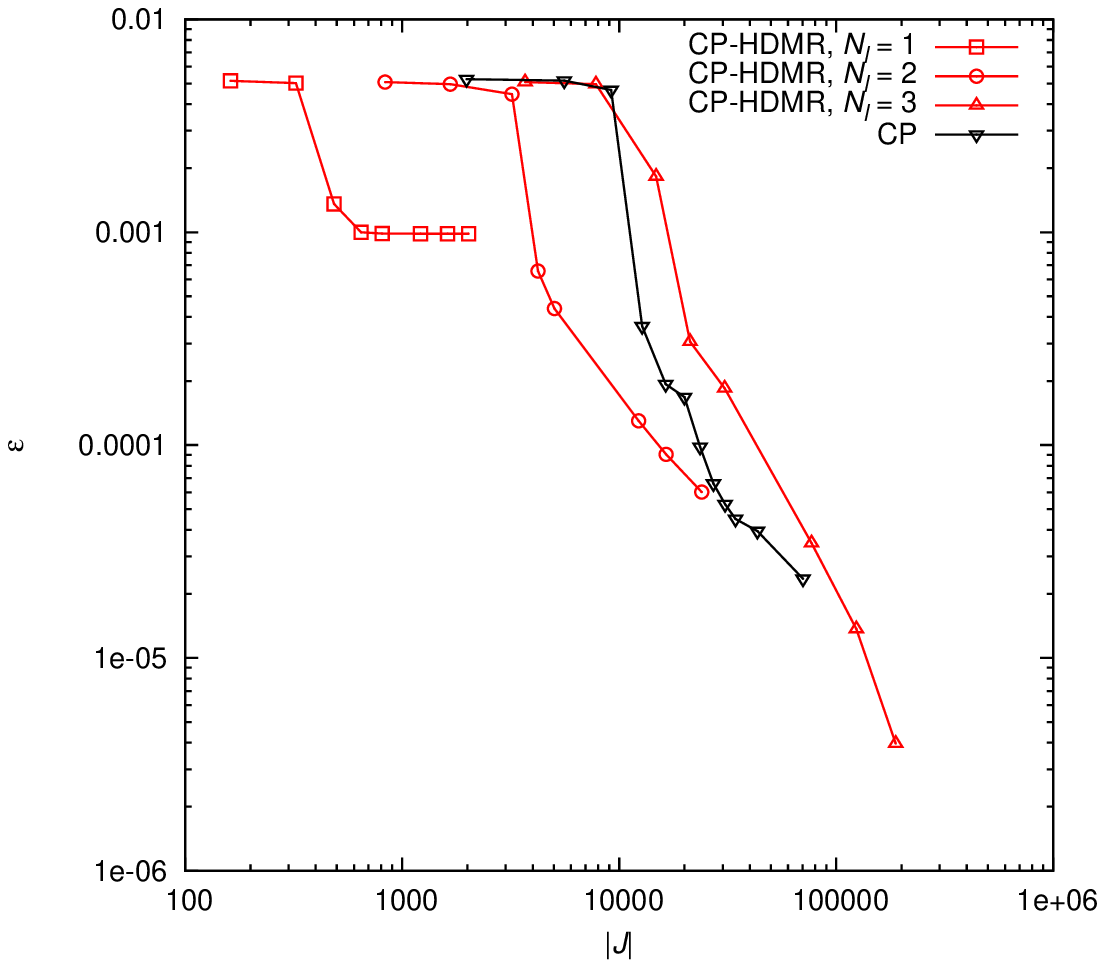}
  \includegraphics[angle = -90, width = 0.49\textwidth, draft = false]{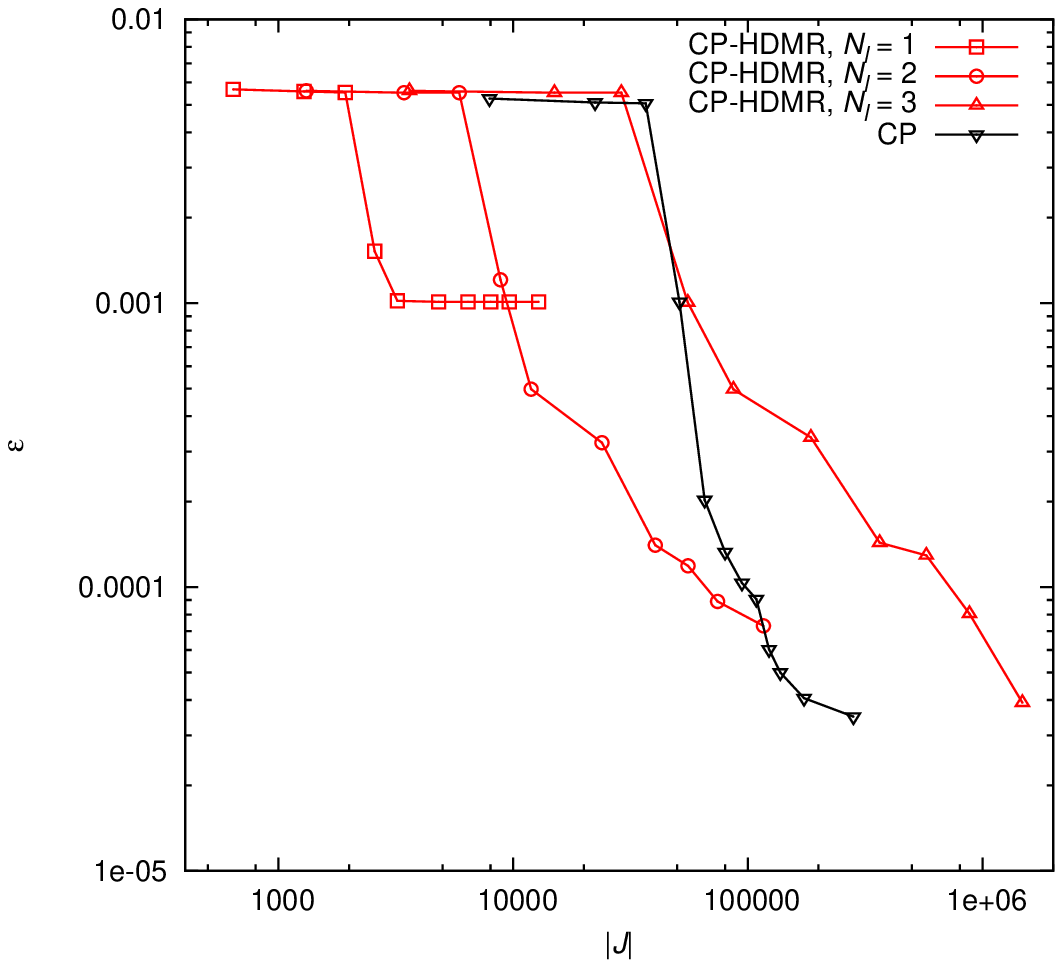}
  \caption{Convergence of the approximation error estimation $\varepsilon$ with the total cardinality $\cardP$ of the representation basis. Approximations of the stochastic modes with the CP-HDMR and CP-like format are compared. $\Nd = 10$ (left) and $\Nd = 40$ (right).}
  \label{L2_card_HDMRvsRK1}
\end{center}
\end{figure}

This CP-HDMR approximation of the stochastic modes is now compared with a CP-like approximation in the form of Eq.~\eqref{CP-like}. To evaluate the CP decomposition, we use an algorithm similar to that in \cite{Nouy_HD_10}. Both decompositions rely on the same approximation basis for the deterministic modes $\left\{\modex_n\right\}$.
We focus on the accuracy of the reconstruction as a function of the cardinality of the whole approximation basis both for $\Nd = 10$ and $\Nd = 40$ when the maximum decomposition rank $\KLrank$ varies, see Fig.~\ref{L2_card_HDMRvsRK1}. The total cardinality increases as more terms are considered in the decomposition series.

The accuracy of the representation is seen to improve as more terms are considered in the series expansion and both the CP and the CP-HDMR formats exhibit an exponential convergence with the total size $\cardP$ of the decomposition.
The first order CP-HDMR quickly reduces the error but plateaus as the functional space is small.
The third order CP-HDMR decomposition is more costly in number of coefficients to evaluate to reach a given error and, in the present settings, should only be considered if high accuracy is needed. 

Unless the targeted accuracy is really high, this motivating example tends to indicate that, for a reasonable required accuracy, a CP-HDMR format involves fewer unknowns than a CP-like decomposition, both for a low $\Nd = 10$- and a moderate $\Nd = 40$-dimensional problem. This is an important point since the number of coefficients which can be evaluated with a reasonable accuracy from experimental data is directly related to the size of the available dataset.
Finally, a CP-HDMR format allows a great flexibility in representing interaction modes $\left\{f_\bgamma\right\}$. In particular, a more parsimonious representation is used in the main text and achieves a similar accuracy with a lower number of terms.
 


\section{Statistics and sensitivity analysis} \label{Sobol_sec}

Once an approximation of a random variable $u\left(\bxi\right)$ is obtained, it is easy to estimate its first statistical moments.
From the HDMR format properties, the estimated mean is simply given by the first term of the decomposition: $\left< u \right>_{L^2\left(\Xi, \, \measxi\right)} \simeq f_\emptyset$.

Thanks to the orthogonality property of the modes $\left\{f_\bgamma\right\}$, the variance $\Var\left(u\right) := \left< \left(u - \left< u \right>_{L^2\left(\Xi, \, \measxi\right)}\right)^2\right>_{L^2\left(\Xi, \, \measxi\right)}$ approximates as the sum of the variance of the individual interaction modes:
\bea
\Var\left(u\right) & \simeq & \sum_{\bgamma \in \basissetfcur \backslash \emptyset}{\Var\left(\fhat_\bgamma\right)}, \nonumber \\
& = & \sum_{\substack{\bgamma \in \basissetfcur \backslash \emptyset \\ |\bgamma| \le \NinterPC}}{\sum_{\balpha, |\balpha| \le \No}{\coef_{\bgamma, \balpha}^2 \, \left\|\psi_\balpha\right\|^2_{L^2\left(\Xi, \, \measxi\right)}}} + \sum_{\substack{\bgamma \in \basissetfcur \backslash \emptyset \\ \NinterPC < |\bgamma| \le \Ninter }}{\sum_{r,r'=1}^{\nr\left(\bgamma\right)}{\prod_{i \in \bgamma}{\sum_{\alpha=1}^{\No}{c_{\bgamma,\alpha}^{r,i} \, \coef_{\bgamma,\alpha}^{r',i} \, \left\|\psi_\alpha\right\|_{L^2\left(\Xi, \, \measxi\right)}^2}}}}, \label{HDMR_var}
\eea
where use was made of the orthogonality of the Hilbertian basis $\left\{\psi_\alpha\right\}$.

Other standard statistical quantities are the sensitivity indices $\left\{S_\bgamma\right\} := \Var\left(\fhat_\bgamma\right) / \Var\left(\uhat\right)$ which essentially represent the relative part of the variance of the QoI due to the interaction of a given set of input random variables only, \cite{Sobol_93,Homma_Saltelli_96}. From Eq.~\eqref{HDMR_var}, it immediately follows that $\sum_{\bgamma \subseteq \left\{1, \ldots, \Nd\right\} \backslash \emptyset}{S_\bgamma} = 1$ and the explicit expression of the sensitivity indices is straightforward to derive from the HDMR format.
In practice, it is often more useful to assess the influence of a given input onto the variance of the QoI with the total sensitivity indices $\left\{S_{T, i}\right\}_{i=1}^{\Nd}$:
\be
S_{T, i} := \frac{\sum_{\bgamma \subseteq \left\{1, \ldots, \Nd\right\} \backslash \emptyset: \, i \in \bgamma}{\Var\left(\fhat_\bgamma\right)}}{\Var\left(\uhat\right)}, \qquad 1 \le i \le \Nd.
\ee

Again, using Eq.~\eqref{HDMR_var}, this quantity is straightforward to estimate once the approximation of $u\left(\bxi\right)$ is available.

\end{document}